\documentclass[letterpaper, reqno,11pt]{article}
\usepackage[margin=1.0in]{geometry}
\usepackage{amsmath,amssymb,amsthm}
\usepackage{graphicx}
\usepackage{wrapfig}
\graphicspath{{images_cropped/} }

\usepackage[percent]{overpic}

\usepackage{color}
\usepackage{comment}

\newtheorem{thm}{Theorem}[section]
\newtheorem{lem}[thm]{Lemma}

\newtheorem{defn}[thm]{Definition}
\newtheorem{conjecture}[thm]{Conjecture}

\newtheorem{cor}[thm]{Corollary}

\newtheorem{rmk}[thm]{Remark}

\newcommand{\RR}{\mathbb{R}}

\newcommand{\PP}{\mathbb{P}}

\newcommand{\EE}{\mathbb{E}}

\newcommand{\pts}{\mathcal{P}}
\newcommand{\lines}{\mathcal{L}}

\newcommand{\cell}{\mathcal{C}}

\begin{document}

 \title{Structure of cell decompositions in Extremal Szemer\'edi-Trotter examples}
\author{Nets Hawk Katz\thanks{Rice University, Houston TX. This work was supported by a grant from the Simons Foundation [SFI-MPS-SIM-00017353, N. H. K.].},\and Olivine Silier\thanks{University of California Berkeley, partially supported by the Jane Street Graduate Research Fellowship}}

\maketitle
\begin{abstract}
The symmetric case of the Szemer\'edi-Trotter theorem says that any configuration of $N$ lines
and $N$ points in the plane has at most $O(N^{4/3})$ incidences. We describe a recipe involving just
$O(N^{1/3})$ parameters which sometimes (that is, for some choices of the parameters)
 produces a configuration of N point and N lines. (Otherwise, we say the recipe fails.) We show that any
near-extremal example for Szemer\'edi Trotter is densely related to a successful instance of the recipe.
We obtain this result by getting structural information on cell decompositions for extremal Szemer\'edi-Trotter examples. We obtain analogous results for unit circles.
\end{abstract}

\section{Introduction}\label{introSection}

If $l$ is a line and $p$ a point in the real plane ${\mathbb R}^2$, we say that
$(l,p)$ is an incidence if $p \in l$.
The most fundamental result in the theory of incidences between points and lines
in the plane is the Szemer\'edi-Trotter theorem \cite{ST} which bounds their number:

\begin{thm}[Szemer\'edi--Trotter 1983]  Let $\lines$ be a set of $n$ lines in the plane
and $\pts$ be a set of $m$ points. Then if $I(\lines,\pts)$ is the set of incidences between
lines of $\lines$ and points of $\pts$, we have the bound
$$|I(\lines,\pts)| \lesssim n^{{2 \over 3}} m^{{2 \over 3}} + n + m.$$
\end{thm}

One thing that is remarkable about the Szemer\'edi-Trotter theorem is that
as far as the exponents are concerned, it is sharp. A number of examples are
known, but we are far from classifying all possible examples. To do so remains
one of the central open problems in incidence geometry of the plane \cite{G}. We restrict
ourselves to the symmetric case where $n=m$, although the question is interesting
whenever $\sqrt{n} < m < n^2$.

\bigskip

{\bf Inverse Szemer\'edi--Trotter problem} Let $\epsilon >0$, $\lines$ be a set of $n$ lines and $\pts$
be a set of $n$ points with 
$$I(\lines, \pts) \geq n^{{4 \over 3}-\epsilon} .$$
What can be said about the structure of $\lines$ and $\pts$ in terms of $n$ and $\epsilon$?

\bigskip

A related question which motivates this study is the unit distance problem.

\bigskip

{\bf Unit distance problem} Let $\pts$ be a set of $n$ points in the plane. Let
$U(\pts)$ be the set of pairs of points in $\pts$ which are at Euclidean distance 1. What upper
bound can one put on $\max_{|\pts|=n}|U(\pts)|$?

\bigskip

Erd\H{o}s introduced the unit distance problem in 1946 \cite{E} and provided the bounds $n^{1 +\frac{1}{\log \log n}} \leq \max_{|\pts|=n}|U(\pts)| \leq n^{3 \over 2}$. He conjectured that the lower bound was close to sharp.

\begin{conjecture}[Erd\H{o}s Unit Distance Conjecture] For all $\epsilon > 0$

\[ |U(\pts)| \leq |\pts|^{1+\epsilon} \]

\end{conjecture}

A recent result by OpenAI provided a counter example to this long standing conjecture with $U(\pts) \gtrsim n^{1+c}$ for a fixed constant $c > 0$ \cite{OAI}. It is currently unclear what the correct conjecture should be.

\bigskip

The best known upper bound for the unit distance problem is $max_{|\pts| = n} |U(\pts)| \leq n^{{4 \over 3}}$. This is not a coincidence [forgive the pun]. Unit distances
are incidences between the points of $\pts$ and the unit circles centered at
those points. Now unit circles are not lines, but they do share some properties in common.
Each unit circle is defined by two parameters and while it is not the case that unit circles
intersect in at most one point, they do intersect in at most two. Essentially every technique
which has been used in a proof of the Szemer\'edi-Trotter theorem can be applied in the case  of unit distances and this is the source of the $n^{{4 \over 3}}$ bound.

\bigskip

A connection between the unit distance problem and the inverse Szemer\'edi-Trotter
problem is that if one had an inverse theorem for unit distances at the exponent
${4 \over 3}$, one could gain a small improvement in the exponent by showing that
the inverse cases don't exist. This is, in fact, a big part of our motivation which is why
we don't mind restricting to the symmetric case in Szemer\'edi-Trotter.

\bigskip

To illustrate the source of the difficulty in obtaining an inverse Szemer\'edi-Trotter
theorem, we describe a simpler, related problem in which the inverse theorem
is fairly straightforward to obtain. We note that the Szemer\'edi-Trotter theorem
uses a great deal more about the structure of the plane than the fact that two
lines intersect at a simple point. If we had restricted ourselves to using only that
fact, we would have obtained this weaker result.

\bigskip

\begin{thm} [Cauchy-Schwarz]   Let $\lines$ be a set of $n$ lines in the plane
and $\pts$ be a set of $m$ points. Then if $I(\lines,\pts)$ is the set of incidences between
lines of $\lines$ and points of $\pts$, we have the bound
$$|I(\lines,\pts)| \lesssim n^{{1 \over 2}} m + n .$$
\end{thm}

\bigskip

The inverse Cauchy-Schwarz problem is to describe all sets of $n$ lines and $n$ points
with $n^{{3 \over 2}}$ incidences. ``It's easy," the reader should exclaim, ``there
are none by Szemer\'edi-Trotter." But we will suspend disbelief and nevertheless try
to describe them despite their nonexistence. What follows is a sketch.

\bigskip

In a configuration of $n$ points and $n$ lines with $n^{{3 \over 2}-}$ incidences,
the typical point is incident to $n^{{1 \over 2}-}$ lines and the typical line is incident
to $n^{{1 \over 2}-}$ points.  We pick an initial point $p_1$. Let $B(p_1)$, the 
``bush" of $p_1$ be the set of points incident to one of the lines incident to $p_1$.
We should have
$$|B(p_1)| \gtrsim n^{1-}.$$
Already, a substantial proportion of the point set $\pts$ belongs to $B(p_1)$ and lies
on the lines going through the point $p_1$. We can go one step further and do this 
twice. We can choose points $p_1$ and $p_2$ so that
$$|B(p_1) \cap B(p_2)| \gtrsim n^{1-}.$$
In other words, a substantial proportion of the point set consists of points lying on a line incident to $p_1$
and a line incident to $p_2$.  After a projective transformation sending $p_1$ and $p_2$
to points at infinity, we get that a substantial portion of the point set lies
on a product set $A \times B$ with each of $|A|$ and $|B|$ of size $n^{{1 \over 2}-}$.
This is not yet an inverse theorem, but it is what we refer to as a \textit{proto-inverse theorem}. Recall that an inverse theorem gives a complete characterization of the solution set to the inverse problem. A proto-inverse theorem on the other hand gives a looser characterization which must include all solutions to the inverse theorem but may also include non-examples.
We have parametrized (a substantial portion of) the point set and whereas {\it a priori}
we needed $O(n)$ parameters to describe the point set, we now need just
$O(n^{{1 \over 2}})$ parameters. This is an important step which has hitherto not been
available in the case of Szemer\'edi-Trotter.

\bigskip

To go from the proto-inverse theorem for Cauchy-Schwarz to an actual inverse theorem
we consider the case of $n^{1-}$ lines having at least $n^{{1 \over 2}-}$ incidences each
with a product set $A \times B$ with each of $A$ and $B$ having size at most 
$n^{{1 \over 2}}$. By rescaling, we can have $0,1 \in A$ with $n^{1-}$ lines having
an incidence with each of $\{0\} \times B$ and $\{1\} \times B$. Thus the lines are 
identified with pairs of points to which they are incident $(0,b_1), (1,b_2)$.  If the same
line is incident to a point of $\{a\} \times B$, we get that $a b_2 + (1-a) b_1  \in B$.
Thus for a typical $a \in A$, the quotient ${1-a \over a}$ has $n^{{3 \over 2}-}$ representations as a member of ${B-B \over B-B}$. This is true for at least 
$n^{{1 \over 2}-}$ choices of $a$. For subsets of the reals, this
phenomenon  is ruled out by the sum-product theorem. In other settings, (finite fields,
the $\delta$-discretized setting) things are bit more delicate because the Szemer\'edi-Trotter
theorem isn't true. Inverse Cauchy-Schwarz, although it didn't go by that name, played
an important role in the development of sum-product theory in those settings. (See 
\cite{KT} and \cite{BKT}.)

\bigskip

In this paper, we obtain the first, to our knowledge, proto-inverse theorem for Szemer\'edi-Trotter and for the unit distance problem at the exponent ${4 \over 3}$.

\begin{thm} There is collection $A$ of $N^{{1 \over 3}}$ parameters and maps  $\lines$
and $\pts$ so that for some values of the parameters $A$,  $(\lines(A),\pts(A))$ is a configuration of at least $N^{1-}$ and at most $N$ lines and points. If $(\lines,\pts)$ is an
extremal configuration of between $N^{1-}$ and $N$ lines and points for Szemer\'edi-Trotter then so is $(\lines \cap \lines(A), \pts \cap \pts(A))$. The analogous result is true for unit circles. \end{thm}

We obtain the parametrization in the theorem (which appears in Theorems \ref{protoinverseST} and \ref{protoinverseSTcircles} in the body of the paper) by
a deep study of the cell decompositions which prove the Szemer\'edi-Trotter theorem.
There is a rather strong analogy to the proto-inverse theorem for Cauchy-Schwarz mentioned
above. In the Cauchy-Schwarz setting most of the points lie on a pair of bushes, or after
a projective transformation, a product set. In the Szemer\'edi-Trotter setting, it is the cell
decomposition which is given by two bushes. We give this as  Theorem \ref{doublebushmixing} and the analogous result for unit circles as Theorem 
\ref{doublebushmixingcircle}. The main idea is that we combine the ideas of cell
decomposition by choosing random lines and of the crossing number inequality. By counting
the crossings inside cells we are able to organize extremal examples using heuristics suggested by random selections.

\section{Extremal examples and cell decompositions}\label{cellSection}

We shall be concerned in this paper with ``extremal examples" for the Szemer\'edi-Trotter theorem.  Our
examples will be (almost) symmetric consisting of  approximately $N$ lines and approximately $N$ points. We will use the notation that the inequality
$$A \lesssim B,$$
between two non-negative quantities $A$ and $B$ will mean that there is a constant $C$ independent of $N$ so that
$$A \leq CB.$$

We would like to allow ourselves losses of small powers of $N$. We choose at the beginning of the paper a
small exponent $\epsilon_0$. Implicitly, at each line of the paper, there will be a different exponent $\epsilon$
depending on the line of the paper, with each $\epsilon$ having the property $\epsilon \lesssim \epsilon_0$.
We will abbreviate $A \lesssim N^{O(\epsilon_0)}$ by $A \lesssim N^{+}$ or $A \lesssim N^{0+}$ .
Similarly we introduce $A \gtrsim N^{-}$ to mean $AN^{O(\epsilon_0)}  \gtrsim 1 $  and $A \sim N^{\pm}$ to mean
$$N^{-} \lesssim A \lesssim N^+.$$
We let exponents add in the natural way.

Our definition of an extremal example for the Szemer\'edi-Trotter theorem will allow $N^+$ errors. This will
be slightly unusual for the study of point-line incidences in the plane.  The reason is that the Szemer\'edi-Trotter theorem is totally sharp in the $\lesssim$ sense.  For that reason, the main tools used in studying the problem have been honed to be sharp in the $\lesssim$ sense. 
 However, there are two reasons we will allow
$N^+$ errors.  The first is that what we're really after are inverse theorems and these will be stronger
and more useful if they apply to  examples that fail to be sharp by $N^+$. 
 The second reason is that we will be studying the properties of probabilistically constructed cell decompositions. While these have been refined to the
$\lesssim$ level, the probabilistic construction for doing that is a bit more sophisticated, and we will be
taking advantage of the ease of use of the simpler one.

\begin{defn}[Extremal configuration]    \label{extreme}  We say that with ${\cal L}$ a collection of at most $N$
and at least $N^{1-}$ lines in the plane
and ${\cal P}$ a collection of at most $N$ and at least $N^{1-}$ points in the plane and with ${\cal I} ({\cal L},{\cal P})$ denoting the set of incidences
(that is pairs $(L,P)$ of a line from ${\cal L}$ and a point from ${\cal P}$ with the point on the line), we say that
$({\cal L},{\cal P})$ is an {\bf extremal configuration} if
$$|{\cal I} ({\cal L},{\cal P})|  \gtrsim N^{{4 \over 3}-}.$$
\end{defn}

Sometimes, we wish to restrict our attention to only a large subset of the incidences of an
extremal configuration. We introduce the notion of an extremal partial configuration.

\begin{defn}[Extremal partial configuration]    \label{partial}  We say that with ${\cal L}$ a collection of at most $N$
and at least $N^{1-}$ lines in the plane
and ${\cal P}$ a collection of at most $N$ and at least $N^{1-}$ points in the plane and with ${\cal I} ({\cal L},{\cal P})$ denoting the set of incidences
(that is pairs $(L,P)$ of a line from ${\cal L}$ and a point from ${\cal P}$ with the point on the line)  and with
$${\cal J}({\cal L},{\cal P}) \subset {\cal I}({\cal L},{\cal P})$$
we say that
$({\cal L},{\cal P}, {\cal J}(\lines,\pts))$ is an {\bf extremal partial configuration} if
$$|{\cal J} ({\cal L},{\cal P})|  \gtrsim N^{{4 \over 3}-}.$$
\end{defn}

Any pair $({\cal L},{\cal P})$  with $\lines$ consisting of lines in the plane and with $N^{1-} \lesssim |\lines| \leq N$
 and with $\pts$ consisting of points in the plane with   $N^{1-} \lesssim |\pts| \leq N$ will be called a {\bf configuration}.

\begin{defn}[Cell decomposition, line weighted]  \label{celll}   Given a configuration $({\cal L},{\cal P})$, we
say that a partition of ${\cal P}$ into $r^2$ disjoint subsets (called cells) $C_1,\dots , C_{r^2}$ is a \textbf{line weighted cell decomposition} if no line in $L \in {\cal L}$ is incident to points in $\gtrsim r$ cells and no cell has 
$\gtrsim {N^{1+} \over r}$ lines of ${\cal L}$ incident to any of its points.   A decomposition having all these properties except the bound on the number of cells a line can be incident to points in will be called a \textbf{provisionally
line weighted cell decomposition}.
\end{defn}

\begin{defn}[Cell decomposition, point weighted]  \label{cellp}   Given a configuration $({\cal L},{\cal P})$, we
say that a partition of ${\cal P}$ into $r^2$ disjoint subsets (called cells) $C_1,\dots , C_{r^2}$ is a \textbf{point weighted cell decomposition} if no line in $L \in {\cal L}$ is incident to points in $\gtrsim r$ cells and no cell contains
$\gtrsim {N^{1+} \over r^2}$ points of ${\cal P}$.
\end{defn}
 
Next we show that any extremal configuration with a line weighted cell decomposition into  approximately $N^{{2/3}}$
parts can be refined into an extremal configuration with a point weighted cell decomposition using the same partition. This is true because cells with too many points do not produce enough incidences per point because of the bound on the number of lines and the Szemer\'edi-Trotter theorem. We simply remove the points of those cells.

\begin{thm} \label{refine1} Let $({\cal L},{\cal P})$ be an extremal configuration with a line weighted cell decomposition $C_1,\dots C_{r^2}$ with $N^{{1 \over 3}-} \lesssim r \lesssim N^{{1 \over 3}}$. Then there is a subset ${\cal P}^{\prime}$
 of ${\cal P}$ with $|{\cal P}^{\prime}| \gtrsim N^{1-}$ and $({\cal L},{\cal P}^{\prime})$ an extremal configuration
so that the nonempty elements of the list $C_1 \cap {\cal P}^{\prime}, \dots, C_{r^2} \cap {\cal P}^{\prime}$ is a point weighted cell decomposition.
\end{thm}
 
\begin{proof}  For the remainder of this proof, we fix the value of $\epsilon$ which corresponds to the current
line of the paper.   We have $|I(\lines,\pts)| \gtrsim N^{{4 \over 3}-\epsilon}$,  we have
$N^{1-\epsilon} \lesssim    |\pts|,|\lines| \leq N$ and we have   
$N^{{1 \over 3} - \epsilon} \lesssim  r \lesssim N^{ {1 \over 3}}$. We divide the cells into two classes
${\cal C}_{big}$ and ${\cal C}_{notsobig}$ where $C_j$ is placed into ${\cal C}_{big}$ if $|C_j|> N^{ {1 \over 3} + 10\epsilon}$ and into ${\cal C}_{notsobig}$ otherwise.  It suffices to take 
$$\pts_{\epsilon}^{\prime} = \bigcup_{C \in {\cal C}_{notsobig}} C,$$ 
and show that 
$$| \pts_{\epsilon}^{\prime} | \gtrsim  N^{1-20\epsilon},$$
and
$$|I(\lines,\pts_{\epsilon}^{\prime} )| \gtrsim N^{ {4 \over 3} - \epsilon}.$$
[This is because at the end of the proof, we can reset the value of $\epsilon$ to $20\epsilon$.]
We calculate

$$|I(\lines,\pts)|=\sum_{C \in {\cal C}_{big}}| I(\lines,C) |+ \sum_{C \in {\cal C}_{notsobig}}| I(\lines,C)|.$$
To bound the first term, we apply the Szemer\'edi-Trotter theorem to each big cell 
using the fact that there are at most $N^{{2 \over 3}+\epsilon}$ lines going through each cell (this is by definition of a line weighted cell decomposition) 
obtaining
$$\sum_{C \in {\cal C}_{big}} |I(\lines,C)| \lesssim \sum_{C \in {\cal C}_{big}} N^{{4 \over 9} +   {2 \over 3}\epsilon}
|C|^{{2 \over 3}}$$
$$ \lesssim \sum_{C \in {\cal C}_{big}} N^{{1 \over 3} - {8 \over 3} \epsilon} |C| \lesssim
N^{{4 \over 3 }- {8 \over 3}\epsilon}.$$
Here the penultimate inequality uses that each $|C|$ is at least $N^{{1 \over 3} + 10 \epsilon}$ and the last inequality uses that $|\pts| \lesssim N$.

Now we know that
$$N^{{4 \over 3} -\epsilon}  \lesssim I(\lines,\pts_{\epsilon}^{\prime} ) ,$$
and we need only show that this implies a good lower bound on $|\pts_{\epsilon}^{\prime} |$. But this follows
immediately from the Szemer\'edi-Trotter theorem and the extremality of the example.

\end{proof}

Next, we will show that for any extremal configuration together with a point weighted cell decomposition with
$N^{{1 \over 3}-} \lesssim r \lesssim N^{{1 \over 3}}$ there is a refinement of the set of lines preserving extremality so that each line
is incident to points in $\gtrsim N^{{1 \over 3}-}$ cells. This is a direct application of
the Cauchy-Schwarz inequality. The lines we remove don't account for many incidences.

\begin{thm}   \label{refine2}
Let $(\lines,\pts)$ be an extremal configuration. Let $C_1, \dots , C_{r^2}$ be a point weighted cell decomposition with $r \sim N^{{1 \over 3} \pm}$. Then there is a refinement $\lines^{\prime} \subset \lines$ so that
$|I(\lines^{\prime},\pts)| \gtrsim N^{{4 \over 3}-}$ and every $L \in \lines^{\prime}$ is incident to points
in $\gtrsim N^{{1 \over 3}-}$ cells.
\end{thm}

\begin{proof}  For the remainder of the proof, we fix the value of $\epsilon$ corresponding to this line,
with $I(\lines,\pts) \gtrsim N^{{4 \over 3}-\epsilon}$.

Consider $\lines_{\epsilon}$, the set of lines intersecting fewer than $r^{1-20\epsilon}$ cells. It 
suffices to show that $|I(\lines_{\epsilon},\pts)|$ is considerably smaller than $N^{{4 \over 3}-\epsilon}$. For each line
$L$, let $C_L$ denote the set of cells in which $L$ is incident to a point. For $L$ a line and $P$ a point, we let
$I_{L,P}$ be the indicator function of incidence,  namely $I_{L,P}=1$ if $P$ is incident to $L$ and $0$ otherwise.

We calculate

\begin{flalign*}
|I(\lines_{\epsilon},\pts)|  &=\sum_{L \in \lines_{\epsilon}} \sum_{C \in C_L} \sum_{P \in C} I_{LP} &\\
                                        &\lesssim N^{{1 \over 2}} r^{{1-20 \epsilon  \over 2}} (
                                              \sum_{L \in \lines_{\epsilon}} \sum_{C \in C_L} 
                                                 (\sum_{P \in C} I_{LP})^2  )^{{1 \over 2}}      &\\
                                      &\lesssim  N^{{1 \over 2}} r^{{1-20 \epsilon  \over 2}} (
                                           I(\lines_{\epsilon},\pts) + \sum_C  \sum_{P_1 \in C} \sum_{P_2 \in C , P_2 \neq P_1} \sum_{L}
                                          I_{LP_1} I_{LP_2}  )^{{1 \over 2}}  &\\
                                   &\lesssim N^{{1 \over 2}} r^{{1-20 \epsilon  \over 2}}  (N^{{4 \over 3}+2 \epsilon})^{{1 \over 2}}   &\\
                                     &\lesssim  N^{{4 \over 3} -{4 \epsilon \over 3}+O(\epsilon^2)}
\end{flalign*}
\end{proof}

The goal of the next theorem is to say that for an extremal example having a point weighted cell decomposition
with $r$ just on the low side of $N^{{1 \over 3}}$,  most of the incidences come from cells having around
${N \over r^2}$ points, and around $ {N \over r}$  lines making just a few incidences with these points, but at least two.
As a corollary, we will obtain a kind of inverse theorem for the lines having a few incidences with the points of a
cell that will prove useful later. (The idea is that for any set of points, the number of lines intersecting two of them is controlled
by the square of the number of points.)

\begin{thm}   \label{refine3}
Let $(\lines,\pts)$ be an extremal configuration. Specifically let $|I(\lines,\pts)|=N^{{4 \over 3}-\epsilon}$ with
$\epsilon$ fixed.
 Let $C_1,\dots , C_{r^2}$ be a point-weighted cell decomposition
for $(\lines,\pts)$ 
with  $N^{{1 \over 3} - 5\epsilon} \leq   r \leq { |I(\lines,\pts)| \over 100 |\lines|}$.  Then there is a set of
incidences  $J(\lines,\pts) \subset I(\lines,\pts)$ so that $|J(\lines,\pts)| \gtrsim  N^{{4 \over 3}-\epsilon}$, but
for every line $L$ and cell $C$ for which there is $P \in C$ with $(L,P) \in J(\lines,\pts)$,   we have that
$$2 \leq | I( \{L\},\pts \cap C)| \lesssim N^+.$$
\end{thm}

\begin{proof}  The way this proof will work is that we will remove from $I(\lines,\pts)$ all incidences that would
violate the conditions for $J(\lines,\pts)$ and observe that we have removed less than half of the 
set $I(\lines,\pts)$.

For any line $L$ and cell $C$ for which there is a unique point $P$ with $(L,P)$ an incidence,
we remove these incidences and we have removed at most $r |\lines|$ incidences since each line has incidences
with at most $r$ cells. Thus $2 \leq | I( \{L\},\pts \cap C)|$.

Next we show that $| I( \{L\},\pts \cap C)| \leq N^+$. We remove incidences from lines $L$ in a cell $C$ if $L$ takes $\gtrsim N^{\epsilon}$ incidences in $C$. To show this did not affect the total incidence count, we must show that rich lines cannot
contribute most of the incidences inside a cell. The number of lines passing through $k$ of the points in a cell is at most
${|\pts \cap C|^2 \over k^3}$, together contributing at most ${|\pts \cap C|^2 \over k^2}$ incidences (this follows from Szemer\'edi--Trotter). So we removed $\lesssim \#(C) {|\pts \cap C|^2 \over N^{2\epsilon}} \lesssim N^{{4 \over 3}-4 \epsilon}$ many incidences. (Here we used that $\sqrt{\#(C)} = r \leq N^{{1 \over 3}-\epsilon}$ and that the cell decomposition is point weighted so $|\pts \cap C| \lesssim N^{{1 \over 3}}$.)

\end{proof}

\begin{cor}  \label{structuredcells}  Let $(\lines,\pts)$ be an extremal configuration. Let $C_1,\dots , C_{r^2}$ be a point-weighted cell decomposition
for $(\lines,\pts)$
with $r \geq  N^{{1 \over 3}-}$ but $r \leq { |I(\lines,\pts)| \over 100 |\lines|}$.   Then there is a set
${\cal C}$ of $\gtrsim r^{2-}$ cells so that for each $C \in {\cal C}$, there is a set of lines $\lines_C$ with
$$|\lines_C|  \gtrsim |C|^{2-},$$
and with each $L \in \lines_C$ incident to at least $2$ but $\lesssim N^+$ points in $C$.  
Each set $\lines_C$ has density $\gtrsim N^-$ in the set of lines intersecting two points in $C$.

\end{cor}

\section{The Probabilistic method and cell decompositions} \label{probSection}

There are several different methods known for producing cell decompositions for proving the Szemer\'edi-Trotter theorem. The most modern is polynomial partitioning. In that method, the boundaries of cells are given
by the zero set of a polynomial.  If a cell decomposition is thus obtained, it is naturally point-weighted. This is called the cellular case. The alternative is that the points all lie in the zero set of a fairly low-degree
polynomial. This is called the structured case. In this case we obtain the
Szemer\'edi-Trotter theorem by bounding the intersection of a curve of
bounded degree and a line by the curve's degree. However, in the structured case we don't have a cell decomposition of the set of points. The method of polynomial partition was first introduced by Larry Guth and the first author in resolving the Erd\"os distinct
distances problem in the plane.  \cite{GK}

An older approach is to define a cell decomposition by randomly selecting lines from the configuration's line set. That approach most naturally produces a cell decomposition that is line-weighted. 
This seems to have been first developed in the seminal paper of Clarkson {\it et. al.}   \cite{CEGSW} as a simplification and improvement of a deterministic construction found in
the original paper of Szemer\'edi and Trotter \cite{ST}. This always produces a cell decomposition, at least if we started with a configuration containing no overly rich lines.

Our present aim is to use cell decompositions in order to learn about
the properties of extremal configurations.  For this purpose, the deterministic nature of polynomial 
partitioning is unhelpful.  The cells are chosen by the Borsuk-Ulam theorem in a somewhat mysterious way. 
Not many choices are available.  We will work with the probabilistic method where almost every selection
of lines yields an acceptable cell decomposition. The fact that an extremal configuration behaves much
 the same for each selection of lines seems to yield a lot of information on extremal configurations. Our original plan for how to carry out our arguments used this observation
heavily. But for technical reasons, it turns out to be beneficial to use a different
classic approach to the Szemer\'edi-Trotter theorem, the crossing number inequality.
Principally, we will use this inequality to say that the lines in a typical point weighted cell
of an extremal configuration behave in the way that you would expect from a random
selection of lines.

We largely follow Terry Tao's blogpost \cite{T} on probabilistic constructions of cell decompositions.  

\begin{defn}
Let $\lines$ be a set
of $N$ lines and $0 \leq r \leq N$. For each line $L \in \lines$ let \[X_r^L=\begin{cases} \{L\} \text{ with probability } {r \over N}\\ \emptyset \text{ otherwise } \end{cases}\]  Let $X_r = \cup_{L \in \lines} X_r^L$. So $X_r^L$ is a random variable valued in subsets of $\lines$ and is defined as a sum of independent i.i.d. Bernoulli random variables $X_r^L$.
\end{defn}

We use one probabilistic calculation repeatedly:

\begin{lem} \label{consecutivelines} Let ${\cal S} \subset \lines$ be any subset of ${C N \log N \over r}$ lines with
$C$ a sufficiently large constant. Then $P(X_r \cap S = \emptyset) \leq {1 \over N^{100}}$.
\end{lem}

\begin{proof}
By independence of the $X_r^L$, the probability that no line in ${\cal S}$ is in $X_r$ is exactly
$$(1-{r \over N})^{{C N \log N \over r}}.$$

By the product limit of the exponential function, as $\frac{N}{r} \to \infty$ the previous expression tends to $N^{-C}$.

\end{proof}

We will also need Chernoff bounds so now is a good time to review them.

\begin{lem} \label{Chernoff}  Let $Y_1, \dots , Y_M$ be independent Bernoulli variables
equal to $1$ with probability $p$ and zero otherwise. Let $pM> M^{\delta}$ for some
$\delta>0$, and let
$$Y=Y_1+Y_2+ \dots Y_M.$$ Then $P(|Y-pM|>pM/10)\lesssim \frac{100}{pM}$. Furthermore $P(Y>10 pM)\lesssim e^{-8pM}$. \end{lem}

\begin{proof}  Observe that since the $Y_j$'s are independent and identically distributed,
we have
$$\EE(e^Y)=\EE(e^{Y_1})^M=(1+p(e-1))^M \sim  e^{p(e-1)M}.$$

In the last $\sim$ we used the product limit of the exponential function. Now we apply Markov's inequality:

$$P(|Y-pM|>pM/10) =  P(|Y-pM|^2>(pM/10)^2) \leq 100\frac{\EE(|Y-pM|^2)}{(pM)^2}. $$

Notice the numerator of the right hand side is just the variance. Recall that independent random variables are uncorrelated so the variance of the sum is equal to the sum of the variance.

$$P(|Y-pM|>pM/10) \leq 100\frac{M \EE(|Y_1-p|^2)}{(pM)^2} \leq 100\frac{M (p(1-p)^2+(1-p)p^2)}{(pM)^2} \leq \frac{100}{pM}.$$

To derive the second bound we again apply Markov's inequality:

$$P(Y>10 pM) = P(e^Y>e^{10 pM}) \leq \frac{\EE(e^Y)}{e^{10pM}}. $$

 recall $\EE(e^Y) \sim  e^{p(e-1)M}$ and $e-1<2$ so $P(Y>10 pM) \leq e^{-8pM}$. Note we were able to obtain a much stronger decay bound for the right tail than for the left tail.

\end{proof}

We are now ready to find a fixed selection $\lines_r \subset \lines$ which will satisfy good properties for partitioning $\lines$ which match the average behavior. We will show that $X_r$ satisfies the intersection of all the good partitioning properties with positive probability. Therefore there exists a specific event which we will call $\lines_r$ that satisfies all the good partitioning properties. First let us explain the three partitioning properties.

\begin{enumerate}
\item Cardinality: We require $\lines_r$ to contain about the average number of expected points: $\frac{r}{2} \leq |\lines_r| \leq 2r $. 

\item For each line $l \in \lines$, we establish an ordering on $\lines \backslash \{l\}$ based on the position of their intersection with $l$. (It might be at infinity.) The ordering is ill-defined when multiple lines are concurrent at a point of $l$, but we order concurrent lines
arbitrarily. 

\begin{defn}\label{defn:partitionsboundarysegments}We say $\lines_r$ \textbf{partitions intersections between lines and boundary segments} if the following two properties hold. For each choice of distinct $l, l' \in \lines$ at least one of the ${C N \log N \over r}$ consecutive lines following $l^{\prime}$ in the order induced by $l$ is in $\lines_r$. Furthermore for each choice of distinct $l, l' \in \lines$ at most $10 \log N$ of the ${\log (N) N \over r}$ consecutive lines following $l^{\prime}$ in the order induced by $l$ is in $\lines_r$.
\end{defn}

If we use $\lines_r$ as boundary lines to partition the plane into disjoint connected components, the first property guarantees that any boundary segment from the partition is intersected by $\lesssim {C N \log N \over r}$ lines from $\lines$. The second property guarantees that at most $10 \log N$ consecutive boundary segments may be sparse. By combining both properties we find that the union of $10 \log N$ consecutive boundary segments must intersect $\sim {N \over r}$ lines from $\lines$.

\item We choose a direction different from that of all lines in $\lines$ which we will refer
to as vertical. At each point $p$ of intersection of two of the lines of $\lines$, the vertical line at $p$ induces an order on the lines
of $\lines$. Again the ordering is ill-defined when multiple lines are concurrent at a point of the vertical line, but we order concurrent lines
arbitrarily. 

\begin{defn}\label{defn:partitionsverticalsegments}We say $\lines_r$ \textbf{partitions intersections between lines and vertical segments} if for each choice of distinct $l, l' \in \lines$ which intersect at a point $p$ at least one of the ${C N \log N \over r}$ consecutive lines above $p$ is in $\lines_r$ and at least one of the ${C N \log N \over r}$ consecutive lines below $p$ is in $\lines_r$.
\end{defn}

\end{enumerate}

\begin{lem}\label{lem:boundary_lines}
Let $\lines$ be a set of $N$ lines and $N^{\delta} \leq r \leq N$. Then there exists a subset $\lines_r \subset \lines$ such that: 

\begin{enumerate}
\item Cardinality: $\frac{r}{2} \leq |\lines_r| \leq 2r $. 

\item $\lines_r$ partitions line-boundary segment intersections. See Definition \ref{defn:partitionsboundarysegments}.

\item $\lines_r$ partitions line-vertical segment intersections. See Definition \ref{defn:partitionsverticalsegments}.

\end{enumerate}

\end{lem}

\begin{proof}
\begin{itemize}

\item Cardinality: We apply Lemma \ref{Chernoff} with $M=|\lines|=N$, $Y_m=1$ if $X_r^{L_m}= \{L_m\}$ and 0 otherwise where $L_m \in \lines$ are all distinct, $p=\frac{r}{N}$, $Y=\# X_r$ and we find that $P(|X_r - r| > r/10) \leq \frac{100}{r}$. So $X_r$ fails to satisfy the cardinality property with probability $\leq \frac{100}{r} \lesssim N^{-\delta}$.

\item Partition line-boundary segment intersections: We first prove that at least one of the next ${C N \log N \over r}$ consecutive lines will be in $X_r$ with high probability. For each choice of distinct $l, l' \in \lines$ we will apply Lemma \ref{consecutivelines} with $\mathcal{S}$ equal the set of ${C N \log N \over r}$ consecutive lines following $l'$ in the order induced by $l$. There are $N^2$ choices of $l, l'$ so the probability that $X_r$ fails to satisfy this property is $\leq N^{-98}$. 

Next we prove that with high probability that among the next ${\log(N)N\over r}$ consecutive lines at most $10 \log(N)$ will be in $X_r$. This follows from the second bound of Lemma \ref{Chernoff} with $M={\log(N)N\over r}$, $p = \frac{r}{N}$, and $Y_m=1$ if $X_r^{L_m}= \{L_m\}$ and 0 otherwise where $L_m \in \lines$ is the $m^{th}$ line in the ordering. We find that the probability that $X_r$ does not satisfy this property is $N^2 P(Y > 10 \log N) \lesssim N^2 e^{-8\log(N)}=N^{-6}$. (The factor of $N^2$ in front accounts for $\leq N^2$ choices of $(l,l')$.)

\item Partition line-vertical segment intersections: For each choice of distinct $l, l' \in \lines$ we will apply Lemma \ref{consecutivelines} with $\mathcal{S}$ equal the set of ${C N \log N \over r}$ consecutive lines above, resp. below the point $l \cap l'$ in the vertical order. There are $N^2$ choices of $l, l'$ so the probability that $X_r$ fails to satisfy this property is $\leq N^{-98}$.

\end{itemize}

$X_r$ fails at least one of these properties with probability $\lesssim N^{-\delta} <<1$ so there exists an event $\lines_r$ which satisfies all three of these properties.
\end{proof}

\begin{defn} \label{funneldecomp} A \textbf{funnel decomposition} is obtained from a cell decomposition by breaking each cell into trapezoids by taking a vertical line segment from each vertex of the cell until the point where it intersects some edge.
(These trapezoids are called funnels in \cite{CEGSW} where this construction was invented.)
\end{defn}

\begin{lem} \label{probconstruction} With $r \sim  N^{{1 \over 3} \pm}$,  the funnel decomposition applied to a set $\lines_r$ from Lemma \ref{lem:boundary_lines} produces a provisionally  line-weighted cell decomposition.
\end{lem}

For our purposes, both the above construction and the slightly refined one of Matou\v{s}ek \cite{M} seem unsatisfactory because of the
introduction of edges which are not on lines of $\lines$. We would like to be able to recognize the lines entering a cell as lines that
intersect a fixed selected line of intersection $L$ between consecutive points of intersection with other selected lines.
For this reason, we find cell decompositions that come from just selecting $r$ random lines and making no further adjustments most natural.
This is problematic because such a decomposition is no longer line weighted due to cells with too many edges. We will deal with
this by bounding the number of points that can be contained in such cells and refining the point set $\pts$ to only include
points in cells with $\lesssim N^+$ edges.

The main ingredient in our bound will be the theorem of Clarkson {\it et. al.} \cite{CEGSW} on the complexity of cell decompositions
given by general families of lines.

\begin{thm} \label{cellcomplexity}  Let $\lines$ be a set of $r$ lines. It divides projective space into $O(r^2)$ cells. Let 
${\cal C}$ be any subcollection of $m$ of these cells. Then the total number of edges of cells in ${\cal C}$ is 
$O(r^{{2 \over 3}} m^{{2 \over 3}} + r)$.   \end{thm}

\begin{cor} \label{numberofbigcells}  Let $\lines_r$ be any set of $r$ lines and let ${\cal C}$ be the set of cells which they
define having $>s$ and $\leq 2s$ edges. Then if $s \leq r^{{1 \over 2}}$, then the number of cells is
$$\#{\cal C} \lesssim {r^2 \over s^3},$$
and if $s \geq r^{{1 \over 2}}$
$$\#{\cal C} \lesssim {r \over s}.$$
\end{cor}

Having now obtained a bound on the number of cells with a certain number of edges, we now control the number of rich points in such
a cell. Let $(\lines,\pts)$ be an extremal configuration in which each point $p$ of $\pts$ is
at least $N^{{1 \over 3}-}$ rich.

\begin{lem} \label{pointsinacell} Let $K$ be a cell coming from a selection of $N^{{1 \over 3}}$ lines $\lines_{N^{{1 \over 3}}}$ from Lemma \ref{lem:boundary_lines}. Suppose
$K$ has $s$ sides. Then $K$ contains at most $sN^{{1 \over 3}+}$ points of $\pts$. \end{lem}

\begin{proof} Rotate the plane such that no line is pointing in the new vertical direction. Take the funnel decomposition of $K$. By Lemma \ref{lem:boundary_lines} the number of lines from $\lines$ that intersect a given vertical segment in the funnel decomposition is $\lesssim N^{{2 \over 3}+}$ and then the number of lines from $\lines$ that intersect a given edge in the boundary of $K$ is $\lesssim N^{{2 \over 3}+}$ so each funnel has $\lesssim N^{{2 \over 3}+}$ lines entering it. Suppose there are $P$ points of $\pts$ in the
funnel $F$. Then there are at least $P N^{{1 \over 3}-}$ incidences in $F$. The Szemer\'edi-Trotter theorem guarantees that
$P \leq N^{{1 \over 3}+}$. Thus the total number of points in $K$ is at most $(s-2) N^{{1 \over 3}+}$ points of $\pts$ which
was to be shown.  \end{proof}

To extract structure from our configurations we will chose to throw out undesirable points and lines keeping only those that enjoy the desirable properties. 

\begin{defn}
    Given an extremal partial configuration $(\lines, \pts, \cal{J(\lines, \pts)})$ we say $(\lines^{\prime}, \pts^{\prime})\subset(\lines, \pts)$ is a \textbf{refinement} if $|\lines| \sim |\lines^{\prime}|$ and $|\pts| \sim |\pts|^{\prime}$ and $|\cal{J(\lines,\pts)}| \sim |\cal{J} (\lines^{\prime},\pts^{\prime})|$.
\end{defn}

Finally, we combine Corollary \ref{numberofbigcells} with Lemma \ref{pointsinacell} to bound the number of points of $\pts$ contained
in cells with between $s$ and $2s$ sides. We call this set of cells $\mathcal{C}$. If $s \leq r^{{1 \over 2}}$ with $r=N^{{1 \over 3}}$, $|{p \in \cup_\mathcal{C} C}| \leq \#\mathcal{C} s N^{{1 \over 3}+} \leq {r^2 \over s^3} sN^{{1 \over 3}+} \leq {N^{1+ } \over s^2}$.  If $s \geq r^{{1 \over 2}}$, we obtain the bound
$|{p \in \cup_\mathcal{C} C}| \leq {r \over s} sN^{{1 \over 3}+} \leq N^{{2 \over 3}+}$. As long as $s$ is much bigger than $N^+$, we do not capture a significant number of points. We conclude the following theorem.

\begin{thm} \label{nicerefinement} Let $(\lines,\pts)$ be an extremal configuration with each point of $\pts$ being at least $N^{{1 \over 3}-}$ rich. Then for each set of lines $\lines_{N^{{1 \over 3}}}$ from Lemma \ref{lem:boundary_lines} there is a refinement $\pts^{\prime} \subset
\pts$ with $|\pts^{\prime}|  \geq {1 \over 2} |\pts|$ so that no point of $\pts^{\prime}$ is contained in a cell with more than $N^{+}$
sides, and in light of Lemma \ref{pointsinacell}, each such cell has at most $N^{{1 \over 3}+}$ points of $\pts^{\prime}$. Furthermore each line in the line set can enter at most $N^{{1 \over 3}}$ cells since it can enter one new cell per time it intersects with a boundary line. Thus we have obtained a point weighted decomposition for the extremal configuration $(\lines,\pts^{\prime})$.
\end{thm}

The main power of Theorem \ref{nicerefinement} for us will be that we can use it to deduce strong properties of extremal examples
without reference to any cell decomposition.

Next we're going to use Corollary \ref{structuredcells} to get structuring result for extremal configurations where many lines
have points bounding intervals where approximately $N^{{2 \over 3}}$ lines intersect. Furthermore large subsets of these groups of about $N^{{2 \over 3}}$ lines are structured: they intersect two of a set of not much more than 
$N^{{1 \over 3}}$ points. A key ingredient in proving this will be the standard crossing number inequality which we state here.

\begin{lem} \label{crossingnumber} \cite{S} Let $G(V,E)$ be a planar graph with $v$ vertices and $e$ edges. Suppose that $e \geq 10v$. Then the number of crossings
between edges is $\gtrsim {e^3 \over v^2}$. \end{lem}


    

We make a definition of a structured set of lines.

\begin{figure}[!h]
    \centering
    \includegraphics[width=\textwidth]{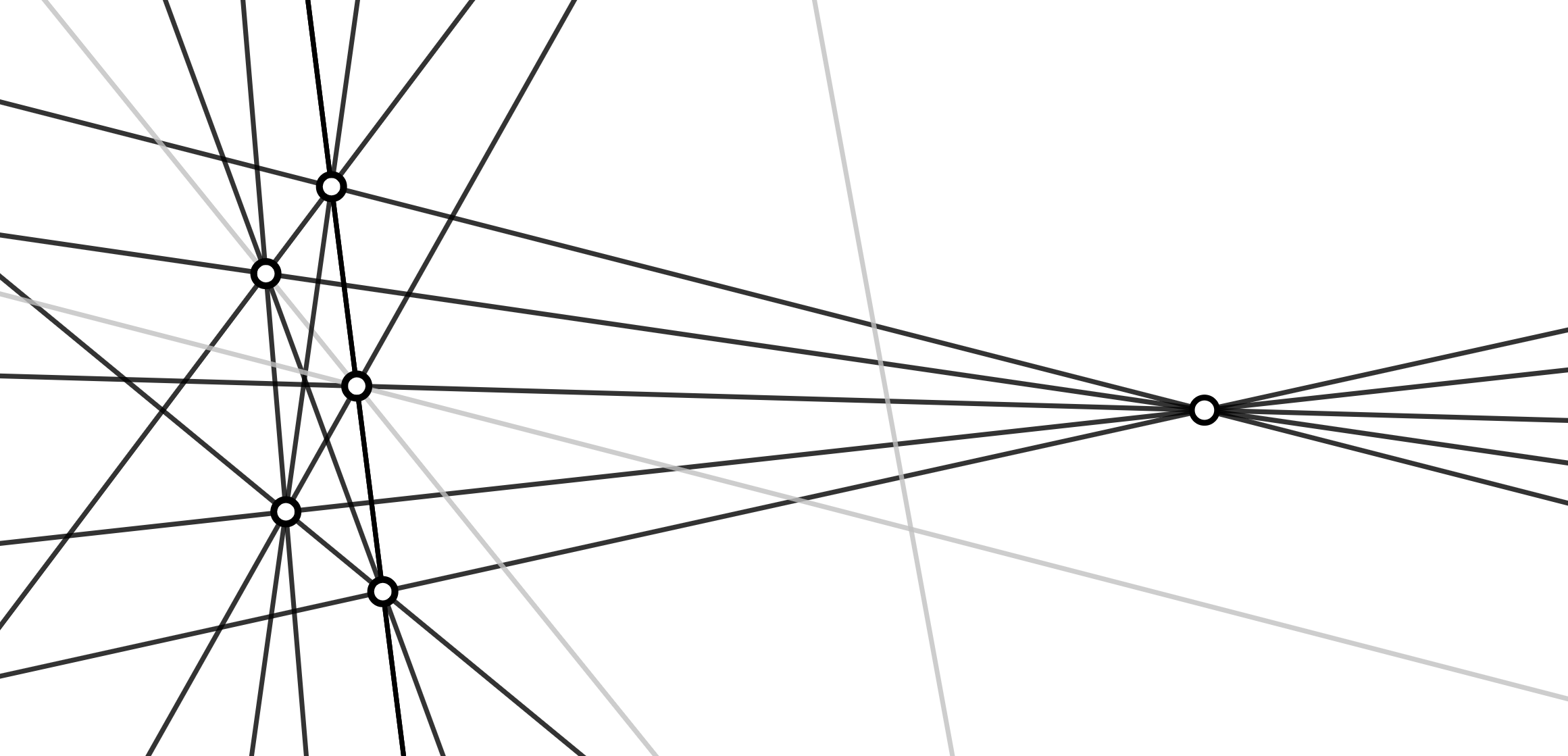}
    \caption{Six structuring points structure the set of black lines. A generic pair of structuring points must define a line in the structured set but not all the pairs must.}
    \label{fig:structured_lines}
\end{figure}

\begin{defn}  We say that a set $\lines_1$ of at least $N^{{2 \over 3}-}$ lines is \textbf{structured} if there is
a set $\pts_1$ of at most $N^{{1 \over 3}+}$ points so that each line of $\lines_1$ is incident to at
least two points of $\pts_1$. We call this set of points \textbf{structuring}. See Figure \ref{fig:structured_lines}. \end{defn}

Note that since $\lesssim N^{\frac{2}{3}+}$ lines go through at least two among $\lesssim N^{\frac{1}{3}+}$ points, the structuring points essentially define the structured lines. Now we're ready to state our structuring theorem.

\begin{figure}[!h]
    \centering
    \includegraphics[width=\textwidth]{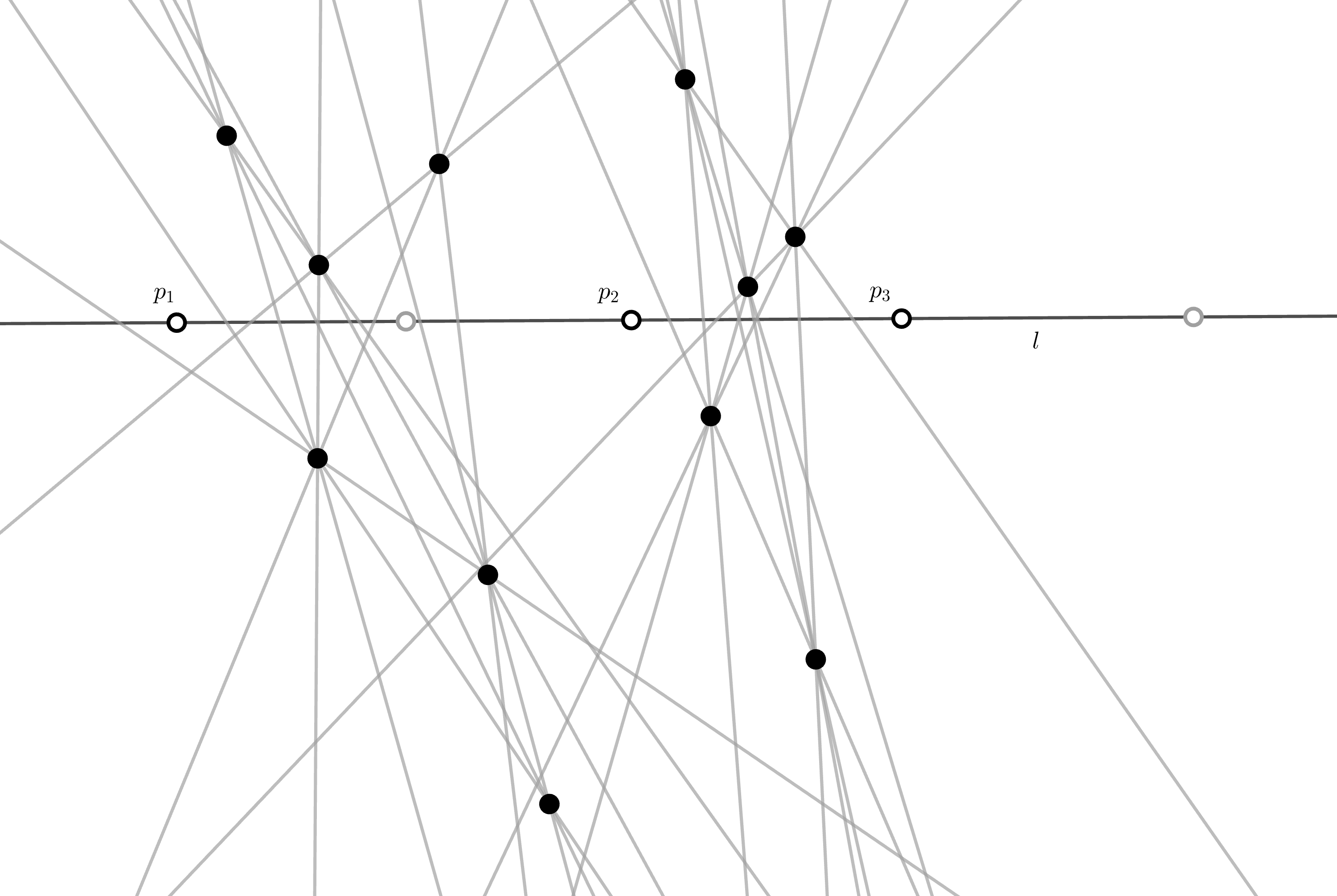}
    \caption{An organizing line $l \in \lines'$ is shown with two intervals bounded respectively by points $p_1,p_2$ and $p_2, p_3$. Each of these intervals is crossed by a structured set of (gray) lines. The two gray points on $l$ are elements of $\pts$ which are not included in our refined construction as the incidence between $l$ and the gray points are not included in $\cal{J}'$.}
    \label{fig:organizing_line}
\end{figure}

\begin{thm} \label{verynicespacing}  Let $(\lines,\pts, \cal{J})$ be an extremal partial configuration. Then there is a refinement $(\lines^{\prime},\pts,\cal{J}')$ so that for each line $l \in \lines^{\prime}$ there are points $p_1, \dots
, p_M$ of $\pts$ with $(l,p_j) \in \cal{J}'$ where $M \gtrsim N^{{1 \over 3}-}$ and the $p_j$'s are ordered by their position on $l$ and so that for each consecutive pair of points $p_j,p_{j+1}$, there is a structured set of lines $\lines_j \subset \lines$ so that each $l^{\prime}$ in $\lines_j$ intersects $l$ in the open interval bounded by the points $p_j$ and $p_{j+1}$. We say the lines in $\lines^{\prime}$ \textbf{organize} $\pts$. See Figure \ref{fig:organizing_line}.\end{thm}

\begin{proof}   Starting with the extremal partial configuration $(\lines,\pts,\cal{J})$ we remove all points with fewer than $N^{{1 \over 3}-}$ incidences in $\cal{J}$
obtaining a refinement $\pts^{\prime} \subset \pts$ so that $(\lines,\pts^{\prime}, \cal{J}^{\prime})$ is still an extremal partial configuration (with $\cal{J}^{\prime}$ the intersection of
$\cal{J}$ with the Cartesian product of $\lines$ and $\pts^{\prime}$)
but 
satisfies the hypotheses of Theorem \ref{nicerefinement}. We apply Theorem \ref{nicerefinement} to obtain a further refinement
$\pts^{\prime \prime} \subset \pts^{\prime}$ so that we have a point weighted cell decomposition ${\cal C} = 
\{ C_1, \dots, C_{r^2} \}$ for the
extremal configuration $(\lines, \pts^{\prime \prime})$  with $r \gtrsim N^{{1 \over 3}-}$ so that $(\lines,\pts^{\prime \prime})$
together with ${\cal C}$ satisfy the hypotheses of Corollary \ref{structuredcells}.  From Corollary \ref{structuredcells} we obtain
a refinement ${\cal C}^{\prime}$ of ${\cal C}$ with $|{\cal C}^{\prime}| \gtrsim N^{{2 \over 3}-}$ and with a subset
$\lines_C$ of  the lines of $\lines$ 
going through any cell $C$ of ${\cal C}^{\prime}$ being a structured set. Hence, any subset $\lines_1 \subset \lines_C$ with
$|\lines_1| \gtrsim N^{{2 \over 3}-}$ is also a structured set.

We now examine a fixed cell $C \in {\cal C}^{\prime}$, the set of points $\pts_C$ consisting of points of $\pts^{\prime \prime}$
which lie inside the cell $C$ and the structured set $\lines_C$ from the previous paragraph.  We define a graph $G(V,E)$
whose vertex set consists of the points of $\pts_C$ and whose edges are the pairs of consecutive points of $\pts_C$ which
are incident to any line of $\lines_C$. Because $\lines_C$ is a structured set (structured by the point set $\pts_C$  because of
the application of Corollary \ref{structuredcells}), there is at least one edge of $G$ for each line of $\lines_C$. Thus,
$v=|V| \lesssim N^{{1 \over 3}+}$ while $e=|E| \gtrsim N^{{2 \over 3}-}$. We conclude from Lemma \ref{crossingnumber} that the
number of crossings for the graph is  $\gtrsim {e^3 \over v^2} \gtrsim N^{{4 \over 3}-}$.  But the crossings of the graph are
nothing other than intersections between lines of $\lines_C$ which occur inside the cell $C$. For each line $l$ in $\lines_C$ which
intersects at least $N^{{2 \over 3}-}$ lines of $\lines_C$, we associate the interval $I$ which is the intersection of the line
with the cell. This interval contains at least one (in fact, two) points of $\pts^{\prime \prime}$ and intersects a structured
subset $\lines_I$ from $\lines_C$. We will count pairs $(l,I)$. Because the lines of $\lines_C$ all must intersect the cell $C$ which
has fewer than $N^{+}$ edges, it must be that $|\lines_C| \lesssim N^{{2 \over 3}+}$ and therefore each cell $C$ of 
${\cal C}^{\prime}$ generates at least $N^{{2 \over 3}-}$ pairs $(l,I)$.  Summing over all the cells of ${\cal C}^{\prime}$,
we obtain at least $N^{{4 \over 3}-}$ many such pairs.  For each line $l$, the different intervals $I$ for which $(l,I)$ are such pairs are
disjoint. Since each interval is crossed by a structured set of lines, there are at most $N^{{1 \over 3}+}$ intervals $I$ for
each line $l$. Thus there must be at least $N^{1-}$ lines $l$ with at least $N^{{1 \over 3}-}$ pairs $(l,I)$. We will
call this set of lines $\lines^{\prime}$.

For each line $l \in \lines^{\prime}$, we order the intervals $I$ for which $(l,I)$ is a pair as $I_1,I_2, \dots ,I_M$.
For each $j =2k-1$ odd, we pick a point $p$ from $\pts$ from the at least two which lie in the interval and call it $p_k$.
For each $j=2k$ odd, we pick the structured set $\lines_k$ which intersects $I_j$.

\end{proof}

We'd now like to take advantage of our result. Theorem \ref{verynicespacing} is a method of associating to each extremal configuration
a refinement which is rather nicely parametrized. We will do this by applying point-line duality to the result of Theorem \ref{verynicespacing}.
The result of the theorem gives us many lines $l$ which are each incident to a set of points $p_1, \dots ,p_M$ with $M \gtrsim N^{{1 \over 3}-}$
so that we have $\gtrsim N^{{2 \over 3}-}$ lines intersecting $l$ between adjacent points which are structured. We use point--line duality, applying Theorem \ref{verynicespacing} to lines and points instead of points and lines.

\begin{figure}[!h]
    \centering
    \includegraphics[width=\textwidth]{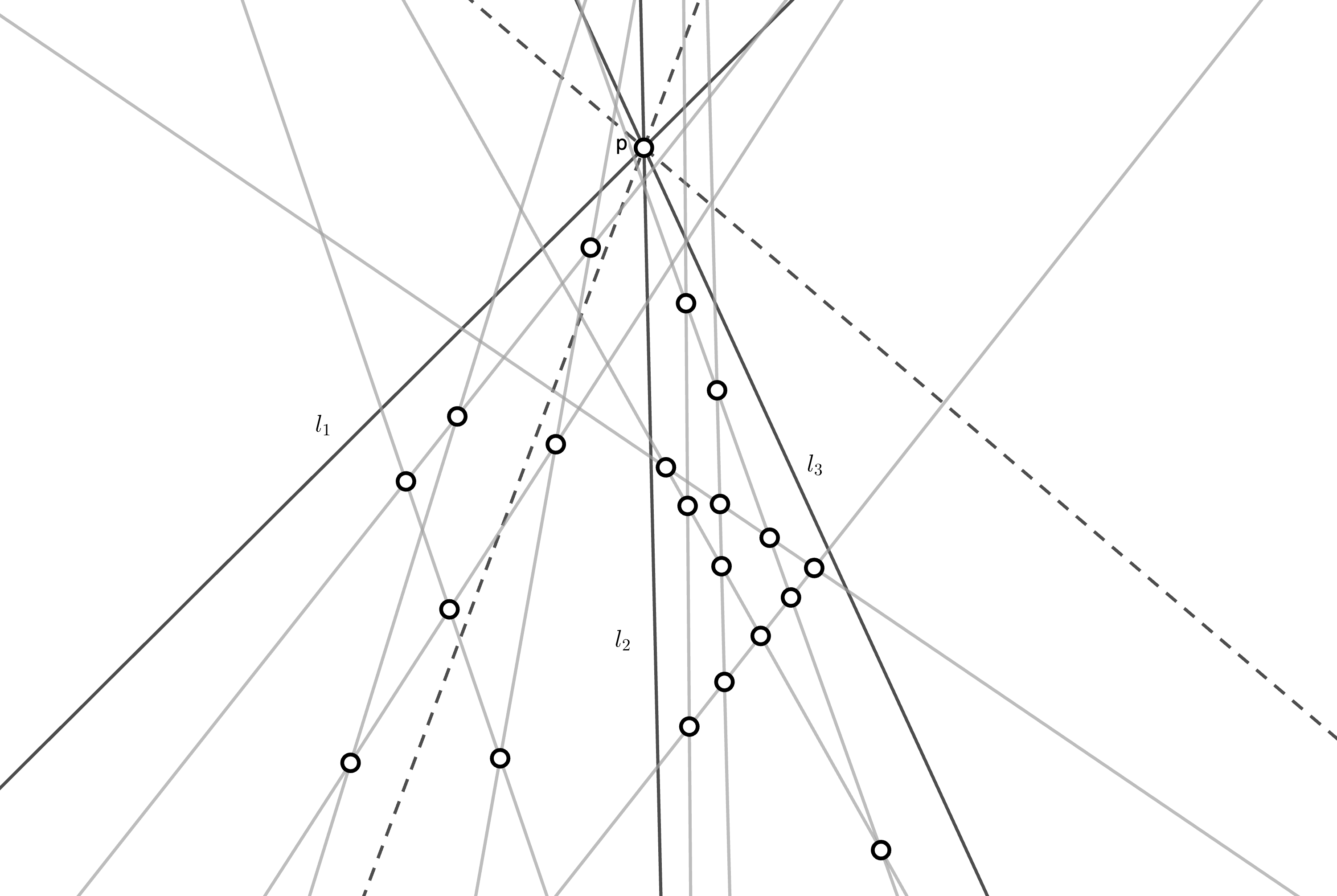}
    \caption{An organizing point $p \in \pts'$ is shown with two sectors $s_1$ and $s_2$ bounded respectively by lines $l_1,l_2$ and $l_2, l_3$. Sectors $s_1$ and $s_2$ each contain a structuring set of (gray) lines. The two dotted black lines through $p$ are elements of $\lines$ which are not included in our refined construction as the incidence between $p$ and the dotted lines are not included in $\cal{J}'$.}
    \label{fig:organizing_point}
\end{figure}

\begin{thm} \label{dualverynicespacing}  Let $(\lines,\pts, \cal{J})$ be an extremal partial configuration. Then there is a refinement $(\lines,\pts^{\prime}, \cal{J}')$ so that for each point $p \in \pts^{\prime}$ there are lines $l_1, \dots
, l_M$ of $\lines$ which are incident to $p$ in $\cal{J}'$ in order of their direction and with $M \gtrsim N^{{1 \over 3}-}$ so that each sector bounded by consecutive pairs of lines $l_j,l_{j+1}$ contains a structured set of $\gtrsim N^{\frac{2}{3}-}$ points $\pts_j \subset \pts$. We say the points in $\pts^{\prime}$ \textbf{organize} $\lines$. See Figure \ref{fig:organizing_point}. \end{thm}

For any point $p \in \pts^{\prime}$ with lines
$l_1, \dots l_M$ incident to it, there are $\gtrsim N^{{2 \over 3}-}$ points in each of the $\gtrsim N^{{1 \over 3}-}$ sectors for a total of $N^{1-}$ points. We take this set of points as a refinement $\pts^{\prime}$ of our original set of points $\pts$. What is particularly pleasant about this structure is that each of the $\gtrsim N^{{2 \over 3}-}$ points
of $\pts^{\prime}$ between two adjacent lines $l_j$ and $l_{j+1}$ lie on at least two of the $\gtrsim N^{{1 \over 3}-}$ structuring lines.

Structuring lines seem very odd precisely because all of the points on them lie in a particular sector between an $l_j$ and $l_{j+1}$. But this is
not as odd as it seems. We see from the proof of Theorem \ref{verynicespacing}, that the structuring lines for $p$ are dual to the points of cells
through which the line dual to $p$ passes. Every point lies in a cell, and every cell has $N^{2/3 \pm}$ lines going through it so by duality the set of $N^{1/3}$ structuring lines define the $N^{2/3 \pm}$ points in a sector.

We're going to show that for any choice of $p \in \pts^{\prime}$ a typical line $l$ will have incidences in most of the sectors between consecutive lines $l_j$ and $l_{j+1}$. To do this we first need to introduce a refinement of the configuration endowed with a cell decomposition whose boundary lines include the bush through $p$. 

\begin{thm}[bush construction] \label{bushconstruction}
For any extremal configuration $(\lines, \pts)$ there exists a subset $\pts^{\prime}$ of $\gtrsim N^{1-}$ points in $\pts$ which are organizing with $\sim N^{{1 \over 3} \pm}$ sectors and a refined configuration $(\lines^{\prime}, \pts^{\prime})$ such that the $\gtrsim N^{1-}$ lines in $\lines^{\prime}$ organize $\pts^{\prime}$. Also for any $p\in \pts^{\prime}$ the refinement $(\lines^{\prime}, \pts_p)$ where $\pts_p$ are the points in $\pts$ organized by $p$, has a refinement
$(\lines^{\prime}, \pts_p^{\prime})$ which is an extremal
configuration with a point-weighted cell decomposition where each cell is contained in a sector. Moreover, any line $l \in \lines^{\prime}$ which crosses exactly $N^{{2 \over 3}+\alpha}$
lines of $\lines^{\prime}$ within the sector $s$ with $\alpha>k\epsilon$ for $k$ sufficiently large
will not enter more than $N^{\alpha+}$ cells in $s$.
 \end{thm}

\begin{proof}

\textit{Refinement properties:} First we keep only points from $\pts$ which are $\sim N^{{1 \over 3} \pm}$ rich. This does not significantly affect the number of incidences. We then apply Theorem \ref{dualverynicespacing} and obtain the refinement $(\lines, \pts^{\prime})$ of organizing points. Then we apply Theorem \ref{verynicespacing} to $(\lines, \pts^{\prime})$ obtaining the refinement $(\lines^{\prime}, \pts^{\prime})$ where the lines in $\lines^{\prime}$ organize $\pts^{\prime}$. Now we have shown the first claim of the theorem.


\textit{Cell decomposition:} Let $p \in \pts^{\prime}$ and $\pts_p$ be the set of points in $\pts$ organized by $p$. Note there are $\gtrsim N^{1-}$ organized points each $\gtrsim N^{{1 \over 3}-}$ rich so $(\lines^{\prime}, \pts_p)$ is an extremal configuration which we work in for this paragraph. We label the bush of $M \sim N^{{1\over 3}\pm}$ lines intersecting $p$ as $l_1,\dots, l_M$.
By Theorem \ref{nicerefinement} we can find a set of lines $\lines_{N^{{1 \over 3}}} \subset \lines'$ which yields a point weighted cell decomposition for some refinement of the point set $\pts_p' \subset \pts_p$. We add in the bush lines $l_1,\dots, l_M$ as additional boundary lines obtaining our final cell decomposition. We need to check this final cell decomposition is still point weighted. Adding $M$ many boundary lines can only reduce the number of points per cell and increases the number of cells traversed by generic lines by $\leq M \lesssim N^{{1 \over 3}+}$. Thus no line in $\lines'$ is incident to points in $\gtrsim N^{{1 \over 3}}$ cells and no cell contains $\gtrsim N^{{2 \over 3}+}$ points so the final cell decomposition is indeed point weighted and its cells are entirely contained in the sectors.

\textit{Crossings:} To get the claim about lines $l$ with $N^{{2 \over 3}+\alpha}$ crossings we use the second part of the $\lines_{N^{1 \over 3}}$ partitions line-boundary segment intersections property in
Lemma \ref{lem:boundary_lines}. This property says that for any list of $\log(N)N^{{2 \over 3}}$ lines which cross $l$ consecutively, at most $10 \log(N)$ of them are in $\lines_{N^{1 \over 3}}$. So $l$ takes $\gtrsim N^{{2 \over 3}-}$ crossings per $10 \log(N)$ consecutive cells and therefore can only intersect $\lesssim N^{\alpha+}$ many cells in the sector.

 \end{proof}

When trying to show that a typical line in a configuration $(\lines,\pts)$ will have incidences in most of the sectors of an organizing point, the enemy case is lines which take too many points in a given sector. So if we show that lines with too many crossings in a sector do not contribute significantly to the total number of incidences in that sector, then most incidences come from lines taking $\lesssim N^{+}$ points in that sector and we win.

\begin{thm} \label{sectorincidences}
There exist $\gtrsim N^{1-}$ organizing points $p$ in $(\lines, \pts)$ with $\sim N^{{1 \over 3}\pm}$ sectors such that for every sector $\gtrsim N^{1-}$ of its incidences come from lines taking $\lesssim N^{-}$ points in that sector.
\end{thm}

\begin{proof}
We choose a point $p$ to be the center of our bush and use the bush construction $(\lines', \pts_p')$ from Theorem \ref{bushconstruction}. First we toss out sectors that have $\gtrsim N^{{5 \over 3}+2\epsilon}$ line-line crossings. There are at least $N^{{1 \over 3}-\epsilon}$ many sectors that have $\gtrsim N^{{5 \over 3}-\epsilon}$ many line-line because there are $\gtrsim N^{2-\epsilon}$ crossings total (if not we can use the crossing lemma version of Szemer\'edi-Trotter to show that the configuration cannot be sharp). So we kept $\sim N^{{1 \over 3}\pm}$ sectors which each have $\sim N^{{5 \over 3} \pm}$ line-line crossings and each sector contributes $\gtrsim N^{1-\epsilon}$ incidences. So this refinement still yields an extremal configuration $(\lines^{\prime}, \pts^{\prime})$.

By Theorem \ref{refine3} we may chose a subset $J(\lines^{\prime}, \pts^{\prime}) \subset I(\lines^{\prime}, \pts^{\prime})$ such that $|J(\lines^{\prime}, \pts^{\prime})| \gtrsim N^{{4 \over 3}-}$ and every line has $\lesssim N^{+}$ incidences from $J(\lines^{\prime}, \pts^{\prime})$ per cell.

\begin{defn}[Fast lines]    \label{fastlines}  We say that a line is $\alpha-$fast for a sector $s$ if it has $\sim N^{{2 \over 3}+\alpha\pm}$ crossings with lines in $\lines$.
\end{defn}

This is the enemy case. An $\alpha-$fast line crosses $N^{{2 \over 3}+\alpha \pm}$ lines of $\lines$ in $s$. By Theorem \ref{nicerefinement} the number of lines that may enter any given cell is $\lesssim N^{{2 \over 3}+}$. Thus an $\alpha-$fast line takes at most $N^{{2 \over 3}+}$ line-line crossings per cell and therefore must enter at least $N^{\alpha-}$ cells. By Theorem \ref{bushconstruction}
each $\alpha-$fast line enters no more than $N^{\alpha+}$ cells in $s$. Our goal is to show that fast lines do not contribute significantly to the number of incidences in the sector. 

\bigskip

Let $l_1$ be the first boundary line of the sector $s$ and let $l_2$ be the second boundary line. Let $i_{1,j}$, resp. $i_{2,j}$ be the set of $\alpha-$fast lines which intersect $l_1$, resp. $l_2$ in the interval bounded by its $j^{th}$ and $(j+1)^{th}$ intersection points with the cell boundary lines from $\lines_r$ (defined in Lemma \ref{lem:boundary_lines}). Note that if $l \in i_{1,j} \cap i_{2,k}$ then $l$ intersects at least $|j-k|$ many cells because $l$ must intersect at least $j-k$ many cell boundary lines which intersect $l_1$ above $i_{1,j}$ and intersect $l_2$ below $i_{2,k}$ (with the orders reversed if $k>j$). Now define \[I_{1,j}=\cup_{k \in [ jN^{\alpha},  (j+1)N^{\alpha}]} i_{1,k} \text{ for } j \in [1, N^{{1 \over 3}-\alpha}]. \]

Define $I_{2,j}$ similarly. Since $\alpha-$fast lines go through $\sim N^{\alpha \pm}$ many cells, they cross $\lesssim N^{\alpha+}$ many boundary lines so $I_{1,j} \cap I_{2,k} = \emptyset$ if $|j-k| \geq N^+$.

\bigskip

Let $l$ be a structuring line for the sector $s$. We inherited the cell decomposition from Theorem \ref{bushconstruction} which guaranteed that it would be point weighted. We can also refine the points to keep only the ones that are $\sim N^{{1 \over 3}\pm}$ rich. There are $\lesssim N$ points so we can refine the cells to keep only cells with $\sim N^{{1 \over 3}\pm}$ points and still obtain an extremal configuration. $l$ intersects $\gtrsim N^{{1 \over 3}-}$ points in $s$ and has $\lesssim N^{+}$ incidences from $\cal{J}(\lines', \pts')$ per cell so $l$ goes through $\gtrsim N^{{1 \over 3}-}$ cells in $s$. We refine the point set keeping only points in the cells that $l$ has an incidence with. Now all the remaining points in the sector are in $\sim N^{{1 \over 3}\pm}$ rich cells skewered along $l$.

\bigskip

Next we refine the set of boundary line segments in $s$ by removing all the $\alpha-$fast boundary lines for $\alpha > 4\epsilon$. The number of fast lines is $\lesssim N^{1-4\epsilon}$ since there are $\lesssim N^{1+}$ line-line crossings total in the sector and each fast line takes $\gtrsim N^{4\epsilon}$ many crossings. So the event that at least half the boundary lines are fast occurs with probability $\lesssim N^{-4\epsilon N^{{1 \over 3}-}}$. We may exclude this unlikely event when selecting our boundary lines from Lemma \ref{lem:boundary_lines}. Thus $l$ still intersects $\gtrsim N^{{1 \over 3}-}$ boundary lines after our refinement. The remaining boundary lines partition $l$ into $\gtrsim N^{{1 \over 3}-}$ many disjoint segments. Each cell contains at most one connected segment of $l$ because for each pair of segments there is a boundary line that crosses $l$ between them, separating them into different half planes. So we kept $\gtrsim N^{{1 \over 3}-}$ cells (which $l$ goes through) so this is still a valid cell decomposition.

\bigskip

Now lines in $I_{1,j}$ intersect $\lesssim N^{\alpha}$ many cells. This is because all the boundary lines which could intersect lines from $I_{1,j}$ must intersect $l_1$ in the interval spanned by $I_{1,j}$, or intersect $l_2$ in the interval spanned $I_{2,j \pm k}$ where $|k| \leq N^+$. This is because all the remaining boundary lines are slow so for $\alpha > 4 \epsilon$ there are no boundary lines that intersect $l_1$ above $I_{1,j}$ and $l_2$ below $I_{2,j\pm N^+}$ or vice-versa. So for a fixed $j$, there are $\lesssim N^{\alpha}$ many boundary lines which contribute to cells that lines in $I_{1,j}$ could intersect. $l$ intersects each of these boundary lines once so $l$ goes through at most one cell per boundary line, thus the lines in $I_{1,j}$ intersect $\lesssim N^{\alpha}$ many cells, so take incidences with $\lesssim N^{{1 \over 3}+\alpha}$ many points.

\bigskip

Now we apply the Szemer\'edi-Trotter Theorem to count the number of incidences from $\alpha-$fast lines by summing over each $I_{1,j}$. Let $\lines_\alpha$ be the set of $\alpha-$fast lines, and $\pts_s$ be the set of points in the sector $s$.

\[|I(\pts_s, \lines_\alpha)| \lesssim \sum_j |I_{1,j}|^{{2 \over 3}} (N^{{1\over 3}+\alpha})^{2 \over 3} \lesssim |\lines_\alpha|^{2 \over 3} (\sum_j (N^{{1\over 3}+\alpha})^2)^{1 \over 3} \lesssim N^{(1-\alpha +){2 \over 3}} (N^{{1 \over 3} - \alpha} N^{{2 \over 3}+ 2\alpha})^{1 \over 3} \lesssim  N^{1-\alpha/3+}.\]

The first inequality follows from the Szemer\'edi-Trotter applied to each interval, the second step follows from H\"older's inequality and $\sum_j |I_{1,j}| = |\lines_\alpha| $. To deduce the third inequality we note there are $\lesssim N^{{5 \over 3}+}$ line-line crossings in $s$ and each $\alpha-$fast line takes $\sim N^{{2 \over 3}+\alpha}$ crossings so $|\lines_\alpha| \lesssim N^{1-\alpha+}$ and we recall $j \lesssim N^{{1 \over 3} - \alpha}$.

The total number of incidences in the sector is $\gtrsim N^{1-\epsilon}$. We showed that $|I(\pts_s, \lines_\alpha)| \lesssim N^{1-\alpha/3}$ so the total number of incidences from fast lines is $\lesssim N^{1-4\epsilon/3}<<N^{1-\epsilon}$. Thus the slow lines contribute the majority of the incidences.

\end{proof}

It is worth emphasizing that Theorem \ref{sectorincidences} says that most (up to $\epsilon$  loss in the exponent) lines will take $O(1)$ points from each sector they crosses. In other words, points on a line are evenly spaced among sectors of the bush. This is the lynchpin fact in our main result (see Section \ref{protoinverseST}).

\begin{figure}[!h]
    \centering
    \includegraphics[width=\textwidth]{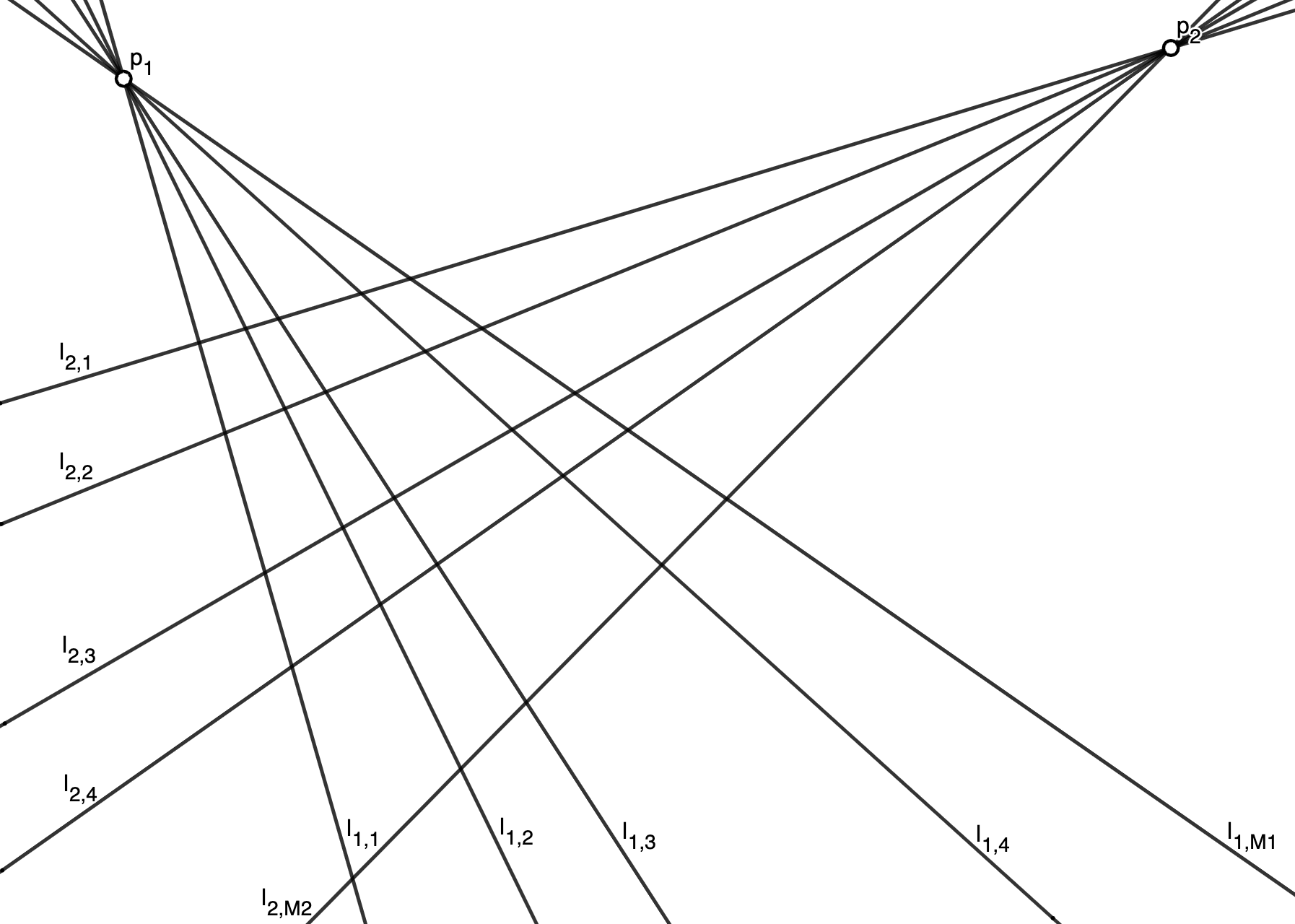}
    \caption{Two organizing points $p_1,p_2$ each have a bush of lines going through them: $l_{1,1},l_{1,2},..., l_{1,M_1}$ and $l_{2,1},l_{2,2},..., l_{2,M_2}$ respectively. These lines generate $M_1M_2$ many point-weighted cells.}
    \label{fig:2bushcelldecomp}
\end{figure}

\begin{thm}[Two Bush Cell Decomposition] \label{doublebushmixing}
Let $(\lines,\pts)$ be an extremal configuration.  Then there are $\gtrsim N^{2-}$ pairs of organizing points $p_1$ and $p_2$ in $\pts$, and two bushes $l_{1,1},\dots, l_{M_1,1}$ incident to $p_1$ and
$l_{1,2}, \dots, l_{M_2,2}$ incident to $p_2$ with $M_1,M_2 \gtrsim N^{{1 \over 3}-}$
and a refinement $\pts^{\prime} \subset \pts$ so that the two bushes break 
$\pts^{\prime}$ into $M_1 M_2$ cells which are point weighted. Moreover each sector $s$
of the bush at $p_1$ with at least $N^{{2 \over 3}-}$ points of $\pts^{\prime}$ has
at least $N^{1-}$ lines of $\lines$ incident to at least 2 points in some cell of the sector. See Figure \ref{fig:2bushcelldecomp}.
\end{thm}

\begin{proof} By Theorem \ref{sectorincidences}, we find a point $p_1$ with bush
$l_{1,1},\dots, l_{M,1}$ and structuring lines holding in total $N^{1-}$ points. We
call this set of points $\pts_1$. Then $(\lines, \pts_1)$ is an extremal configuration.
Also by Theorem \ref{sectorincidences} we find a refinement of the set of incidences
$J(\lines,\pts_1)$ with $|J(\lines,\pts_1)| \gtrsim N^{{4 \over 3}-}$ so that each line of $\lines$
takes only $N^{+}$ incidences of $J(\lines,\pts_1)$ in each sector of the bush at $p_1$.

We restrict to those points in $\pts_1$ which are at least $N^{{1 \over 3}-}$ rich in incidences of $J(\lines,\pts_1)$.  We refer to that set as $\pts_2$. The set $(\lines,\pts_2)$
is an extremal configuration.  We refine the set of lines to $\lines_1$ which take $\lesssim N^+$
incidences in $N^{{1 \over 3}-}$ sectors of the bush at $p_1$. We let $\pts_3$ be
the set of points of $\pts_2$ that are $N^{{1 \over 3}-}$ rich with respect to lines of $\lines_1$ incident to only $N^+$ other points in the same sector. By Theorem \ref{sectorincidences} the pair $(\lines_1,\pts_3)$ is an extremal configuration. Pick an organizing point $p_2$ with bush 
$l_{1,2} \dots l_{M_2,2}$ and structuring lines from $\lines_1$ for each sector of the bush. We call this the second bush.
Let $\pts^{\prime}$ be a set of points that is initially empty. For every sector $s$ from the second bush and for every structuring line in $s$ consider the set of ``good" sectors from the first bush where that structuring line takes $\lesssim N^+$
incidences. Add all the points generated by that structuring line which are in good sectors from the first bush to $\pts^{\prime}$. There are $\lesssim N^{{1 \over 3}+}$ structuring lines in each sector of the second bush, each of which contribute $\lesssim N^+$ points per sector of the first bush so there are $\lesssim N^{{1 \over 3} +}$ points from $\pts^\prime$ in the intersection of a sector from the first bush and a sector from the second bush. Thus the cell decomposition given by the two bushes is point weighted.

Each sector of the first bush with $\gtrsim N^{{2 \over 3}-}$ points has $N^{1-}$ incidences. By Theorem \ref{sectorincidences} we can refine the set of lines incident to this sector to keep only those with $\lesssim N^{+}$ incidences in the sector. By Lemma \ref{refine3} we can apply a final refinement to the line set keeping only those that have at least 2 incidences in some cell in the sector.

\end{proof}

\section{A proto inverse Szemer\'edi-Trotter theorem} \label{protoiSTSection}

A proto inverse Szemer\'edi-Trotter theorem will be a recipe for constructing a configuration
of points and lines which may not terminate or may not yield an extremal configuration but
so that a large portion of every extremal example can be obtained using this recipe.

When considering such a recipe, an important piece of information is how many parameters
one needs to specify to obtain an instantiation of the recipe.

There is a trivial recipe taking $O(N)$ parameters. Namely use $4N$ parameters to 
completely specify a set of $N$ lines, $\lines$ and a set of $N$ points $\pts$. Examine
the set of incidences between these lines and points $I(\lines,\pts)$. If it happens
to be that $|I(\lines,\pts)| \gtrsim N^{{4 \over 3}-}$, then we have constructed
an extremal configuration and in fact, every extremal configuration can be constructed in 
this way. This recipe and its proto inverse Szemer\'edi-Trotter theorem amount to really
just the definition of extremal configuration, and nothing has been gained.

Now, however, we describe a recipe using just $O(N^{{1 \over 3}})$ parameters.
Our recipe will be based on a cell decomposition consisting of a grid of
axis parallel rectangles.  $a_1 < a_2  < \dots < a_{N^{{1 \over 3}}}$ will be real numbers
representing the $x$ coordinates of the grid. $b_1 < b_2  < \dots < b_{N^{{1 \over 3}}}$
will be the $y$ coordinates of the grid. The final ingredients will be a set of lines
$l_{s,1} \dots  l_{s,N^{{1 \over 3}}}$ which will serve as the structuring lines for the
sector of cells between $x$ coordinates $a_1$ and $a_2$.

We declare the recipe to have failed if the lines $l_s$ do not have at least 
$N^{{2 \over 3}-}$ crossings with $x$-coordinate between $a_1$ and $a_2$. Otherwise,
we declare the recipe to have failed unless there are at least $N^{{1 \over 3}-}$ values
of $j$ so that at least $N^{{ 1 \over 3}-}$ and no more than $N^{{1 \over 3}+}$ of the
crossings with $x$-coordinates between $a_1$ and $a_2$ have $y$ coordinates between
$b_j$ and $b_{j+1}$. Otherwise, we define these crossings to be the points of the
cell $[a_1,a_2] \times [b_j,b_{j+1}]$. For each cell, we find all lines which are incident to
two points. We declare this set of lines to be $\lines$. We say that the recipe has failed
unless there are at least $N^{{2 \over 3}-}$ choices of $(j,k)$ so that $N^{{2 \over 3} \pm}$ lines of $\lines$ cross the cell $[a_j,a_{j+1}] \times [b_k, b_{k+1}]$. Otherwise, we
say that the recipe has failed unless for at least $N^{{2 \over 3}-}$ of these choices
$(j,k)$ the lines going through the $(j,k)$th cell are structured. If they are structured,
we refer to the structuring points as the points of the $(j,k)$th cell. And combining all
of these structuring points, we get the set $\pts$ and we declare that the construction
has succeeded. In this case, $(\lines,\pts)$ is an extremal configuration. We denote the
output of this recipe as a function of its inputs as $(\lines(a,b,l_s), \pts(a,b,l_s))$.

First we show that for every extremal example there are parameters for which the recipe will succeed and generate a point set and line set that includes this example with $\gtrsim N^-$ density.

\begin{thm}  \label{protoinverseST}
Let $(\lines,\pts)$ be an extremal configuration. Then there is a choice of the 
$O(N^{{1 \over 3}})$ parameters $(a,b,l_s)$ and a projective transformation $P$ so that
$(P(\lines) \cap \lines(a,b,l_s),P(\pts) \cap \pts(a,b,l_s))$ is an extremal configuration
(of $N^{1-}$ lines $N^{1-}$ points and $N^{{4 \over 3}-}$ incidences.)
\end{thm}

\begin{proof}  This follows from Theorem \ref{doublebushmixing}. We will choose $P$ to be a projective transformation sending the points $p_1$ and
$p_2$ to the points at infinity corresponding to the $y$ direction and $x$ direction respectively. Let $a$ be the list of x coordinates of the vertical lines in the bush $P(p_1)$ and $b$ be the list of y coordinates of the horizontal lines in the bush $P(p_2)$. Let $s$ be
a vertical sector with $\gtrsim N^{{2 \over 3}-}$ points of $\pts$ and $l_s$ be the structuring lines of this sector. These structuring lines structure the $\gtrsim N^{{2 \over 3}-}$ points in $s$. By the last part of Theorem \ref{doublebushmixing} there are $\gtrsim N^{1-}$ lines from the example that are in a set structured by the points in one of the cells of $s$. These are the lines generated by the recipe. By Lemma \ref{refine3} $\gtrsim N^{1-}$ of these lines intersect $\gtrsim N^{{1 \over 3}-}$ cells in which they take at least 2 and $\lesssim N^+$ incidences. So the $N^{{1 \over 3}\pm}$ points in each cell structure the $N^{{2 \over 3}\pm}$ lines through it and so the recipe has succeeded. \end{proof}

Next we will present an equivalent but dual picture of the two bush cell decomposition. Here duality means point line duality.

\begin{defn}
Let $D$ be the map that takes points to lines and lines to points in $\RR \PP^2$ as follows: $D(a,b)=\{y=ax-b\}$ and $D(\{y=ax-b\})=(a,b)$. We call this map point-line duality.
\end{defn}

Note that point-line duality is an involution and preserves incidence structure.

\bigskip

Recall that an organizing line is partitioned into $N^{{1 \over 3}\pm}$ intervals by points it is incident to and a $N^+$ density of the intervals intersect a structured set of $N^{{2 \over 3}\pm}$ lines. These sets of structured lines are dual to the sets of points contained in each sector defined by an organizing point.

\begin{defn} Let $(\lines,\pts)$ be an extremal configuration. Given an organizing line $l$ and its bush of points from $\pts$: $\{p_{1}, \dots, p_{M_1}\}$ on $l$ with $M\gtrsim N^{{1 \over 3}-}$, for each $j \in [1, M]$ we call the subset of $\lines$ which intersects $l_1$ between $p_{j,1}$ and $p_{j+1,1}$ a sector of lines. Note each organizing point defines $\gtrsim N^{{1 \over 3}-}$ many sectors of lines.
\end{defn}

In Theorem \ref{doublebushmixing} we saw that two organizing points and their bush of lines yield a point-weighted cell decomposition where the cells are the intersection of a sector from the first organizing point with a sector from the second organizing point. We similarly define the dual notion of cells of lines:

\begin{defn}Let $(\lines,\pts)$ be an extremal configuration. Given a pair of organizing lines $l_1$ and $l_2$ in $\lines$ and their bush of points from $\pts$: $\{p_{1,1}, \dots, p_{M_1,1}\}$ on $l_1$ and $\{p_{1,2}, \dots, p_{M_2,2}\}$ on $l_2$ with $M_1,M_2 \gtrsim N^{{1 \over 3}-}$ we define $\gtrsim N^{{2 \over 3}-}$ many cells of lines as follows. For each $(j,k) \in [1, M_1] \times [1, M_2]$ the cell of lines $(j,k)$ is the set of lines in $\lines$ that intersects $l_1$ in the interval between $p_{j,1}$ and $p_{j+1,1}$ and intersects $l_2$ in the interval between $p_{k,2}$ and $p_{k+1,2}$.
\end{defn}

\begin{thm} \label{dualprotoinverseST}
Let $(\lines,\pts)$ be an extremal configuration. Then there are $\gtrsim N^{2-}$ pairs of organizing lines $l_1$ and $l_2$ in $\lines$ and a refinement $\lines' \subset \lines$ so that $(l_1, l_2)$ partition $\lines'$ into $\gtrsim N^{{2 \over 3}-}$ many cells of lines each of cardinality $\sim N^{{1 \over 3}\pm}$. Furthermore each of the sectors of lines are structured by $N^{{1 \over 3} \pm}$ structuring points.
\end{thm}

\begin{wrapfigure}[15]{r}{0.4\textwidth}
    \centering
    \includegraphics[width=0.4\textwidth]{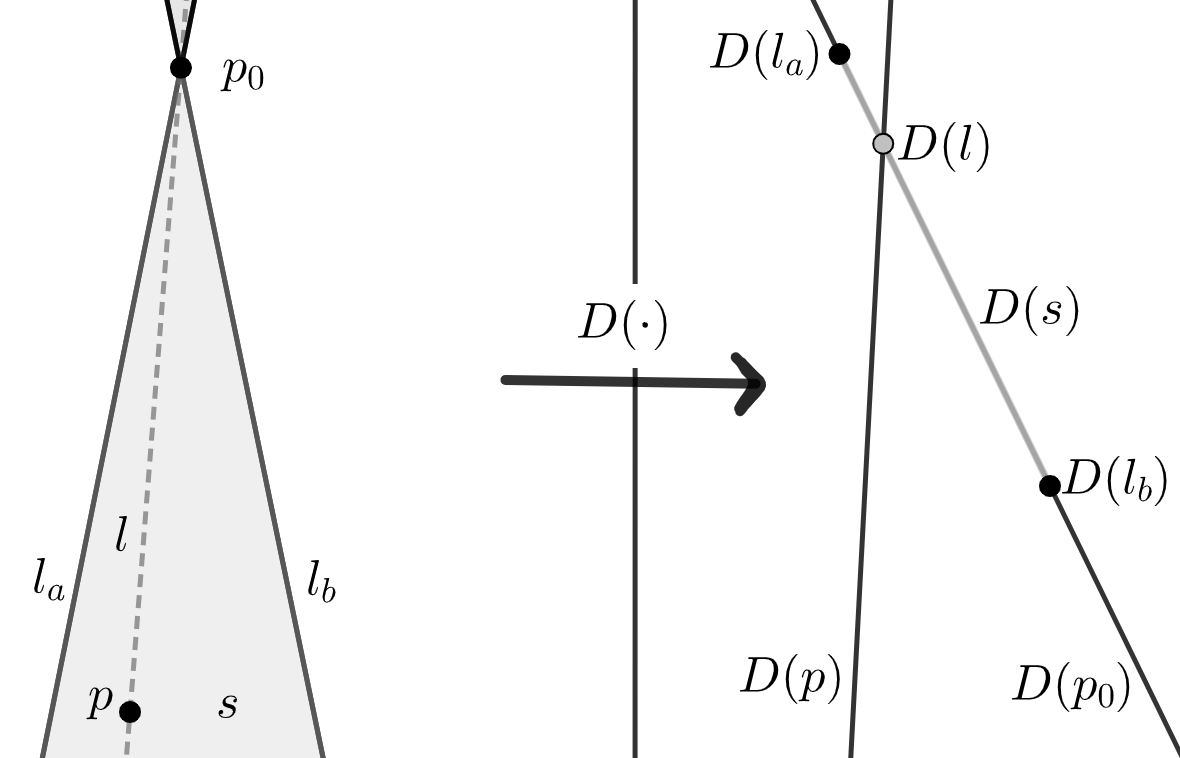}
    \caption{On the left a point $p$ is in the sector $s$ and on the right its dual line $D(p)$ is in the dual line sector.}
    \label{fig:sector_pt_line_duality}
\end{wrapfigure}

\textit{Proof.}  We apply point line duality to the statement of Theorem \ref{doublebushmixing}. First we need to check that the notions of organizing point and organizing line are dual to each other. They are because their properties are entirely defined by point-line incidences which are preserved under point line duality. 

Next we need to check that the notions of sector are dual to each other. In other words we need to check that a point $p$ is in a sector $s$ if and only if $D(p)$ is in the line sector $D(s)$.

For this paragraph see Figure \ref{fig:sector_pt_line_duality}. Let $p_0$ be an organizing point and $p$ be a point in a sector $s$ bounded by the lines $l_a$ and $l_b$. This is equivalent to saying that $p$ lies on a line $l$ through $p_0$ where the angle of $l$ is in the circle arc bounded by the angles of $l_a$ and $l_b$ that corresponds to the interior of $s$. This is equivalent to having $D(p_0)$ an organizing line with a sector $D(s)$ composed of all the lines that intersect $D(p_0)$ in the interval between $D(l_a)$ and $D(l_b)$ such that $D(p)$ intersects the point $D(l) \in D(p_0)$ and the $x$ coordinate of $D(l)$ is in the interval bounded by the $x$ coordinates of $D(l_a)$ and $D(l_b)$. This is equivalent to $D(p)$ is in $D(s)$. (Note $p_0$ has finite $x$ coordinate so $D(p_0)$ has finite slope i.e. is not vertical so the interval spanned by $D(l_a)$ and $D(l_b)$ includes an open interval of $x$ coordinates.)

Since cells are pairwise intersections of sectors the definitions of cells in the two bush cell decomposition are also dual to each other.

The last part of the theorem (sectors are structured by $N^{{1 \over 3}+}$ structuring points) follows from the definition of an organizing line.
\qed

\begin{figure}[!h]
    \centering
    \includegraphics[width=\textwidth]{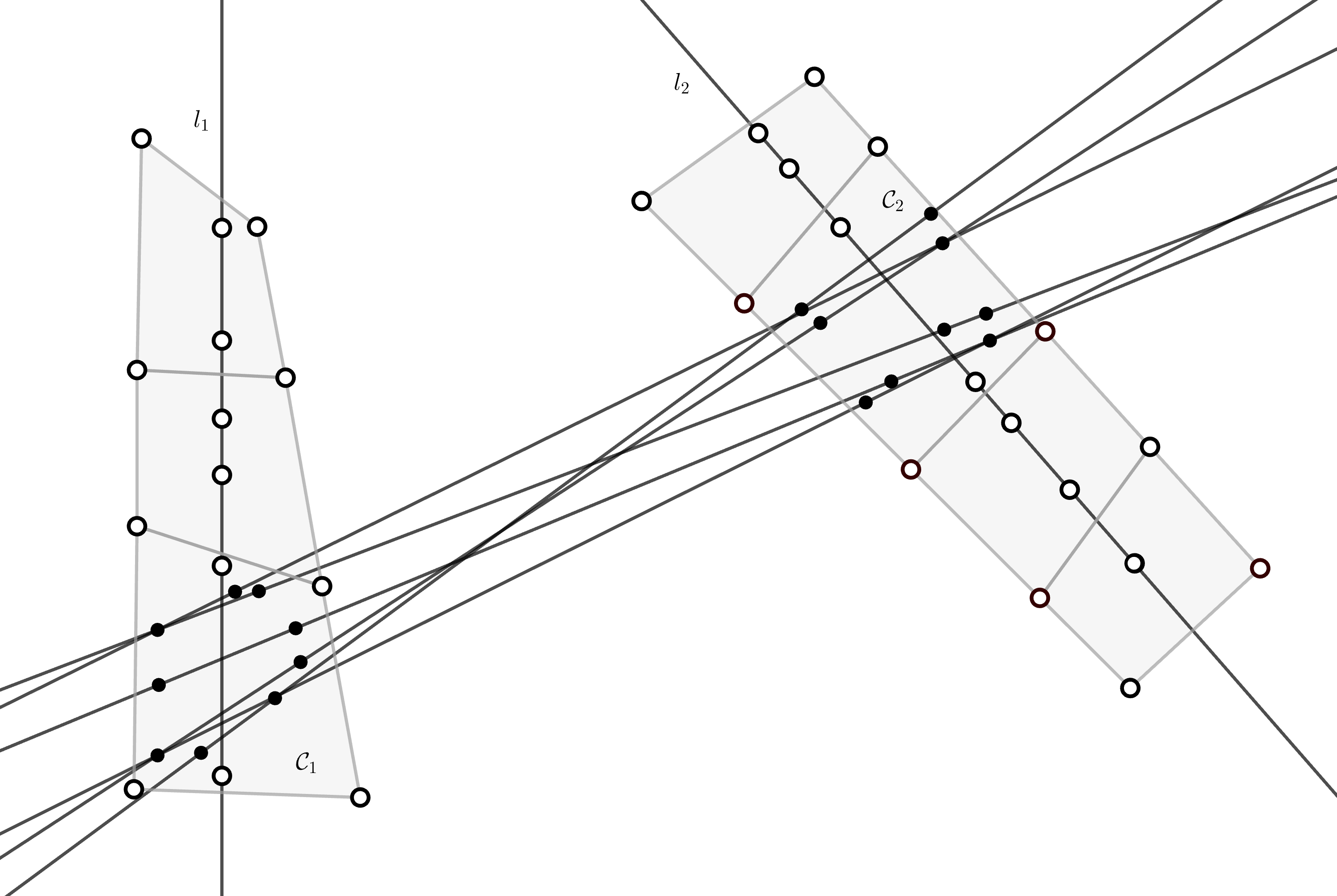}
    \caption{Two organizing lines $l_1,l_2$ go through several cells including cells $\mathcal{C}_{1},\mathcal{C}_{2}$ respectively. These cells are shown to share many lines which each take two points in $\mathcal{C}_{1}$ and in $\mathcal{C}_{2}$.}
    \label{fig:mixing}
\end{figure}

\begin{thm}[Mixing]
\label{mixing}
Let $(\lines,\pts)$ be an extremal configuration with the cell decomposition given by Theorem \ref{protoinverseST}. Then $\gtrsim N^{{4 \over 3}-}$ pairs of cells share $\gtrsim N^{{1 \over 3}-}$ lines which take at least two incidences of $J'$ in each of the two cells. See Figure \ref{fig:mixing}.
\end{thm}

\begin{proof}
    Assume we have the cell decomposition from Theorem \ref{protoinverseST}. An application of Theorem \ref{dualprotoinverseST} gives us that $\gtrsim N^{2-}$ pairs of lines are organizing. 

    Consider a pair of organizing lines $l_1, l_2$. We take a refinement of the line set such that each line has at least two incidences in a cell that $l_1$ has incidences in and at least two incidences in a cell that $l_2$ has incidences in. Note this must yield an extremal refinement for $\gtrsim N^{2-}$ pairs of organizing lines $l_1,l_2$ because each of the $\gtrsim N^{1-}$ organizing lines in $\lines'$ has at least 2 incidences in $\gtrsim N^{{1 \over 3}-}$ cells which each have $\gtrsim N^{{2 \over 3}-}$ lines going through it which each have at least two incidences in that cell. So each organizing line contributes $\gtrsim N^{1-}$ lines (out of a total of $< N$ lines) so $\gtrsim N^{1-}$ other organizing lines must share $\gtrsim N^{1-}$ regular lines that take at least two points in one of each of their cells.

    From our initial application of Theorem \ref{dualprotoinverseST} we know that $\gtrsim N^{{2 \over 3}-}$ pairs of intervals between adjacent points on $l_1$ and adjacent points on $l_2$ share $\gtrsim N^{{1 \over 3}-}$ lines. (Note this is still true after our refinement because we kept $\gtrsim N^{1-}$ lines). Furthermore lines take $\lesssim N^{+}$ incidences in a single cell so $\gtrsim N^{{2 \over 3}-}$ pairs of cells that $l_1$ and $l_2$ take incidences in share $\gtrsim N^{{1 \over 3}-}$ lines out of $\lesssim N^{{2 \over 3}}$ total pairs of cells. 

    Since this result holds for any generic pair of $\gtrsim N^{2-}$ pairs of organizing lines, we conclude that any generic pair of cells must share $\gtrsim N^{{1 \over 3}-}$ lines.
\end{proof}

So we conclude that our two bush cell decomposition can be chosen to be constructed using two sets of $\sim N^{{1 \over 3} \pm}$ parallel lines (projective transformation) which form a grid. $\gtrsim N^{{2 \over 3}-}$ of the rectangles in the grid are cells which contain $\sim N^{{1 \over 3}\pm}$ points and $\sim N^{{2 \over 3}\pm}$ lines which take about 2 points in the cell. Finally we have the mixing property that generic pairs of cells share $\sim N^{{1 \over 3} \pm}$ lines.

\section{Cell decompositions for extremal configurations in the unit distances problem}

The $N^{\frac{4}{3}}$ incidence bound from the Szemer\'edi-Trotter Theorem also holds for points and unit circles in $\RR^2$ and can be proved in the same way \cite{SST}. This bound is not believed to be sharp for unit circles. We hope that progress on the inverse problem for points and unit circles with $\gtrsim N^{\frac{4}{3}-}$ incidences could eventually distinguish between lines and unit circles and lead to an end-point improvement in the unit distance bound.

\begin{thm}[Szemer\'edi--Trotter for unit circles]\label{thm:circle_Sz-trotter} Let $\lines$ be a set of unit circles
and $\pts$ be a set of points both in $\RR^2$. Then if $I(\lines,\pts)$ is the set of incidences between
circles of $\lines$ and points of $\pts$, we have the bound
$$|I(\lines,\pts)| \lesssim |\pts|^{{2 \over 3}}|\lines|^{{2 \over 3}} + |\pts| + |\lines|.$$
\end{thm}

We adapt the argument from Sections \ref{cellSection} to \ref{protoiSTSection} to the incidence problem between points and unit circles centered at the points, which is equivalent to the unit distance problem. We will begin by proving the existence of some point weighted cell decomposition via the probabilistic method e.g. Theorem \ref{nicerefinementcircle}. Next we will define organizing, structuring, and structured sets for point - unit circle arrangements and show that we can find a cell decomposition subordinate to the sectors of an organizing point/ unit circle e.g subsection \ref{subsection:structuring_organizing_bush_circle}. Finally we will show that a pair of organizing points/ unit circles yield a very simple to parametrize and symmetric cell decomposition that we call the two bush cell decomposition e.g. Theorem \ref{doublebushmixingcircle}. We end the section with a discussion of some of the two bush cell decomposition's properties for unit circles.

Section \ref{cellSection} only contains counting arguments that rely on the property that two lines intersect in at most one point. Since two unit circles intersect in at most two points, the counting arguments all still hold in the unit circle case. That said, to tackle some complications due to circle cells lacking convexity, we will slightly modify the chosen number of cells and points per cell. Furthermore the first few probability results in Section \ref{probSection} also rely on counting arguments so also hold for unit circles.

\subsection{Cell Decomposition Preliminaries}
 As before, the unit circle cell decompositions that currently exist in the literature \cite{CEGSW} do not satisfy our purpose because they add additional boundary components to the cells which are not part of our original set of circles. A main ingredient in our proof will be the theorem of Clarkson et. al. \cite{CEGSW} on the complexity of cell decompositions given by general families of unit circles.
 
 \begin{thm} \label{cellcomplexitycircle}  Let $\lines$ be a set of $r$ unit circles. It divides $\RR^2$ into $O(r^2)$ connected components called \textbf{cells}. Let 
${\cal S}$ be any subcollection of $m$ of these cells. Then the total number of edges of cells in ${\cal S}$ is 
$O(r^{{2 \over 3}} m^{{2 \over 3}} \beta(r) + r)$  where $\beta(r)$ is the inverse of the Ackermann function.  \end{thm}

Note the inverse Ackermann function grows notoriously slowly \cite{A} (much slower than $\log(n)$). 


\begin{cor} \label{numberofbigcellscircle}  Let $\lines_r$ be any set of $r$ unit circles and let ${\cal C}$ be the set of cells which they
define having $>s$ and $\leq 2s$ edges. Then the bound from Theorem \ref{cellcomplexitycircle} gives us $|{\cal C}| = O(r^{{2 \over 3}} |{\cal C}|^{{2 \over 3}} \beta(r)/s + r/s)$. If $s \leq r^{{1 \over 2}}/\beta(r)^{3/2}$ then the first term dominates and the number of cells is
$$\#{\cal C} \lesssim \beta(r)^3{r^2 \over s^3},$$
if $s \geq r^{{1 \over 2}}/\beta(r)^{3/2}$ then the second term dominates and we get
$$\#{\cal C} \lesssim {r \over s}.$$
\end{cor}

To extract structure from our configurations we will chose to throw out undesirable points and unit circles keeping only those that enjoy the desirable properties. 

\begin{defn}
    Given an extremal partial unit circle configuration $(\lines, \pts, \cal{J(\lines, \pts)})$ we say 

    \noindent $(\lines^{\prime}, \pts^{\prime},  \cal{J^{\prime}(\lines^{\prime}, \pts^{\prime})} ) \subset (\lines, \pts, \cal{J(\lines, \pts)})$ is a \textbf{refinement} if $|\lines| \sim |\lines^{\prime}|$ and $|\pts| \sim |\pts|^{\prime}$ and $|\cal{J(\lines,\pts)}| \sim |\cal{J}^{\prime}(\lines^{\prime},\pts^{\prime})|$.
\end{defn}

Now will find a set $\lines_{N^{{1 \over 3}}} \subset \lines$ of $\sim N^{{1 \over 3}}$ many unit circles that satisfy good properties for partitioning unit circles in $\lines$ which match the average behavior. 

\bigskip

For each circle $c \in \lines$ we establish an ordering on $\lines$ based on the position of their intersection with $c$. We order the circles that do not intersect $c$ arbitrarily at the end of the ordering. We order the unit circles that do intersect $c$ by choosing a reference point $p \in c$ and a direction (clockwise or counter clockwise) and then ordering them by the order of their first point of intersection with $c$ starting at $p$ and going in the chosen direction. The ordering is ill-defined when multiple circles are concurrent at a point of $c$, but we order concurrent circles
arbitrarily. For each pair of distinct unit circles $(c, c^{\prime}) \in \lines^2$ we consider the order induced by $c$ with a reference point $p \in c \cap c'$. 

\begin{defn} We say $\lines' \subset \lines$ \textbf{partitions intersections between unit circles and boundary arcs} if the following two properties hold. First for each pair of distinct unit circles $(c, c^{\prime}) \in \lines^2$ at least one of the $CN^{{2 \over 3}}\log N$ unit circles following $c'$ in the order induced by $c$ in the clockwise, resp. anticlockwise direction is in $\lines'$.
Second at most $10 \log N$ of the $N^{{2 \over 3}}\log N$ unit circles following $c'$ in the order induced by $c$ in the clockwise, resp. anticlockwise direction is in $\lines'$.
\end{defn}

At each point $p$ of intersection of two unit circles of $\lines$, the vertical line going through $p$ induces an order on the circles
of $\lines$ which is the order of the first point of intersection between a circle in $\lines$ and the vertical line (either in ascending or descending order). We would like to exclude the case that none of the first $C N^{{2 \over 3}} \log N$ circles above $p$ are chosen for $\lines_{N^{{1 \over 3}}}$ and none of the first $C N^{{2 \over 3}} \log N$ circles below $p$ are chosen for $\lines_{N^{{1 \over 3}}}$.

\begin{defn} We say $\lines' \subset \lines$ \textbf{partitions intersections between unit circles and vertical segments} if for each pair of distinct unit circles $(c, c^{\prime}) \in \lines^2$ which intersect at a point $p$ at least one of the next $CN^{{2 \over 3}}\log N$ unit circles above resp. below $p$ is in $\lines'$.
\end{defn}

\begin{lem}\label{lem:boundary_circles} Let $\lines$ be a set of $N$ unit circles. Then there exists a subset $\lines_{N^{{1 \over 3}}}$ such that $\frac{N^{{1 \over 3}}}{2} \leq |\lines_{N^{{1 \over 3}}}| \leq 2N^{{1 \over 3}}$, and such that $\lines_{N^{{1 \over 3}}}$ partitions intersections between unit circles and boundary arcs and such that $\lines_{N^{{1 \over 3}}}$ partitions intersections between unit circles and vertical segments.
\end{lem}

\begin{proof} This proof is the same as that of Lemma \ref{lem:boundary_lines} with $r=N^{{1 \over 3}}$ . The proof translates directly over to circles since the probability Lemmas \ref{consecutivelines} and \ref{Chernoff} do not make any geometry assumptions. In fact this lemma holds for any sets with the same cardinality and ordering definitions.
\end{proof}

\begin{lem} \label{pointsinacellcircle} Let $K$ be a cell coming from a set of unit circles $\lines_{N^{{1 \over 3}}}$ which satisfies the good partitioning properties from Lemma \ref{lem:boundary_circles}. Suppose
$K$ has $s$ sides. Then $K$ contains at most $sN^{{1 \over 3} +}$ many $\gtrsim N^{{1 \over 3}-}-$rich points of $\pts$. \end{lem}

\begin{proof} Following the construction in \cite{CEGSW}, at every intersection of pairs of circles in $\lines_{N^{{1 \over 3}}}$ draw the vertical segment going up until it hits the next circle from $\lines_{N^{{1 \over 3}}}$ and the vertical segment going down until it hits the next circle from $\lines_{N^{{1 \over 3}}}$. The connected components bounded by circles from $\lines_{N^{{1 \over 3}}}$ and vertical segments are called \textbf{funnels}. Each funnel is bounded by at most two circle arcs (top and bottom) and two vertical segments (left and right). $\lines_{N^{{1 \over 3}}}$ has the good partitioning properties guaranteed by Lemma \ref{lem:boundary_circles} so each boundary element of a funnel has at most $C N^{{2 \over 3}} \log(N)$ circles entering it. Suppose there are $P$ points of $\pts$ in the funnel $F$. Then there are at least $P N^{{1 \over 3}-}$ incidences in $F$. The Szemer\'edi-Trotter theorem for unit circles (Theorem \ref{thm:circle_Sz-trotter}) guarantees that
$P \leq N^{{1 \over 3} +}$. 

If $K$ has $s$ sides, then it has at most $s+1$ funnels so the total number of points in $K$ is at most $s N^{{1 \over 3} +}$ points of $\pts$ which
was to be shown.  \end{proof}

Finally, we combine Corollary \ref{numberofbigcellscircle} with Lemma \ref{pointsinacellcircle} to bound the number of $\gtrsim N^{{1 \over 3}-}$ rich points of $\pts$ contained
in cells with between $s$ and $2s$ sides. We call this set of cells $\mathcal{C}$. If $s \leq r^{{1 \over 2}}/ \beta(r)^{{3 \over 2}}$ with $r=N^{{1 \over 3}}$, $|{p \in \cup_\mathcal{C} C}| \leq \#\mathcal{C} s N^{{1 \over 3}+} \leq \beta(r)^3{r^2 \over s^3} sN^{{1 \over 3}+} \leq {N^{1+ } \over s^2}$.  If $s \geq r^{{1 \over 2}}/\beta(r)^{{3 \over 2}}$, we obtain the bound
$|{p \in \cup_\mathcal{C} C}| \leq {r \over s} sN^{{1 \over 3}+} \leq N^{{2 \over 3}+}$. As long as $s$ is much bigger than $N^+$, we do not capture a significant number of points. We conclude the following theorem.

\begin{thm} \label{nicerefinementcircle} Let $(\lines,\pts)$ be an extremal unit circle configuration with each point of $\pts$ being at least $N^{{1 \over 3}-}$ rich. Then for each selection of boundary circles $\lines_{N^{{1 \over 3}}}$ satisfying the good partitioning properties of Lemma \ref{lem:boundary_circles}, there is a refinement $\pts^{\prime} \subset
\pts$ with $|\pts^{\prime}|  \geq {1 \over 2} |\pts|$ so that no point of $\pts^{\prime}$ is contained in a cell with more than $N^{+}$
sides, and in light of Lemma \ref{pointsinacellcircle}, each such cell has at most $N^{{1 \over 3}+}$ points of $\pts^{\prime}$
so that we have obtained a point weighted decomposition for the extremal configuration $(\lines,\pts^{\prime})$.
\end{thm}

\begin{thm} \label{boundedsidescircle}
Let $(\lines,\pts)$ be an extremal unit circle configuration. Specifically let $|I(\lines,\pts)|=N^{{4 \over 3}-\epsilon}$ with
$\epsilon$ fixed.
 Let $C_1,\dots , C_{r^2}$ be a point-weighted cell decomposition
for $(\lines,\pts)$ 
with  $N^{{1 \over 3} - 5\epsilon} \leq   r \leq N^{{1 \over 3} - 4\epsilon}$.  Then there is a set of
incidences  $J(\lines,\pts) \subset I(\lines,\pts)$ so that $|J(\lines,\pts)| \gtrsim  N^{{4 \over 3}-\epsilon}$, but
for every circle $c \in \lines$ and cell $C$ for which there is $P \in C$ with $(c,P) \in J(\lines,\pts)$,   we have that
$$2 \leq | I( \{c\},\pts \cap C)| \lesssim N^+$$ and there exists some other point $P'$ in $C$ such that $P$ and $P'$ are adjacent on $c$ and the circle arc $(P,P')$ is entirely contained in $C$.
\end{thm}

\begin{proof}  The way this proof will work is that we will remove from $I(\lines,\pts)$ all incidences that would
violate the conditions for $J(\lines,\pts)$ and observe that we have removed less than half of the 
set $I(\lines,\pts)$.

For any circle $c$ and any cell $C$ in which $c$ has $\lesssim N^{2\epsilon}$ incidences
we remove these incidences. This procedure removes at most $r |\lines| N^{2\epsilon} \lesssim N^{{4 \over 3}-2\epsilon}$ incidences since each circle has incidences
with at most $r$ cells. By Theorem \ref{nicerefinementcircle} each cell has at most $N^{+} \sim N^{\epsilon}$ sides. Furthermore, every circle intersects a cell side in at most two points. Thus for every circle $c$ and cell $C$, $c$ has at most $N^{\epsilon}$ disjoint connected circle arcs contained in $C$. By the previous refinement, if $c$ has an incidence with a point in $C$ in $J(\lines, \pts)$ then $c$ has at least $N^{2\epsilon}$ incidences in $C$ so by pigeonhole principle, at least two incidences with points $P$ and $P'$ must occur on the same connected circle arc in $C$.

Now we show $| I( \{c\},\pts \cap C)| \lesssim N^+$. We remove incidences between a circle $c$ and points in a cell $C$ if $| I( \{c\},\pts \cap C)| \gtrsim N^{5\epsilon}$ and we must show this did not affect the total number of incidences. The number of circles passing through $k$ of the points in a cell $C$ is at most
${|\pts \cap C|^2 \over k^3}$, together contributing at most ${|\pts \cap C|^2 \over k^2}$ incidences in $C$ (this follows from Szemer\'edi-Trotter for unit circles). Thus we removed at most $(\# C){|\pts \cap C|^2 \over (N^{5\epsilon})^2}$ incidences. Since $(\# C) = r^2 \leq N^{{2 \over 3}-8\epsilon}$ and $|\pts \cap C| \leq N^{{1 \over 3}+}$ by Theorem \ref{nicerefinementcircle}, we kept most of the incidences. Thus $N^{2\epsilon} \lesssim | I( \{L\},C)| \lesssim N^{5\epsilon}$.

\end{proof}

We obtain the following corollary.

\begin{cor}  \label{structuredcellscircle}  Let $(\lines,\pts)$ be an extremal unit circle configuration. Let $C_1,\dots , C_{r^2}$ be a point-weighted cell decomposition
for $(\lines,\pts)$
with $r \sim  N^{{1 \over 3}-}$. Then there is a set
${\cal C}$ of $\gtrsim r^{2-}$ cells so that for each $C \in {\cal C}$, there is a set of circles $\lines_C$ with
$$|\lines_C|  \gtrsim |C|^{2-},$$
and with each $c \in \lines_C$ incident to at least $2$ but $\lesssim N^+$ points in $C$ such that for at least two of these points adjacent on $c$, the circle arc between them is entirely contained in $C$.  
Each set $\lines_C$ has density $\gtrsim N^-$ in the set of unit circles intersecting two points in $C$.

\end{cor}

\subsection{Crossing Numbers and Geometic Preliminaries}

\begin{lem}\label{lem:circle_arc_tube_sqrt}
        A circle arc of length $\sqrt{\delta}<<1$ of a unit circle is contained in a tube of width $\delta$ and length $\sqrt{\delta}$ tangent to the circle (its long axis is perpendicular to the line going through the center of the unit circle and the center of the arc).
\end{lem} 

\begin{wrapfigure}[12]{r}{0.3\textwidth}
    \centering
    \includegraphics[width=0.3\textwidth]{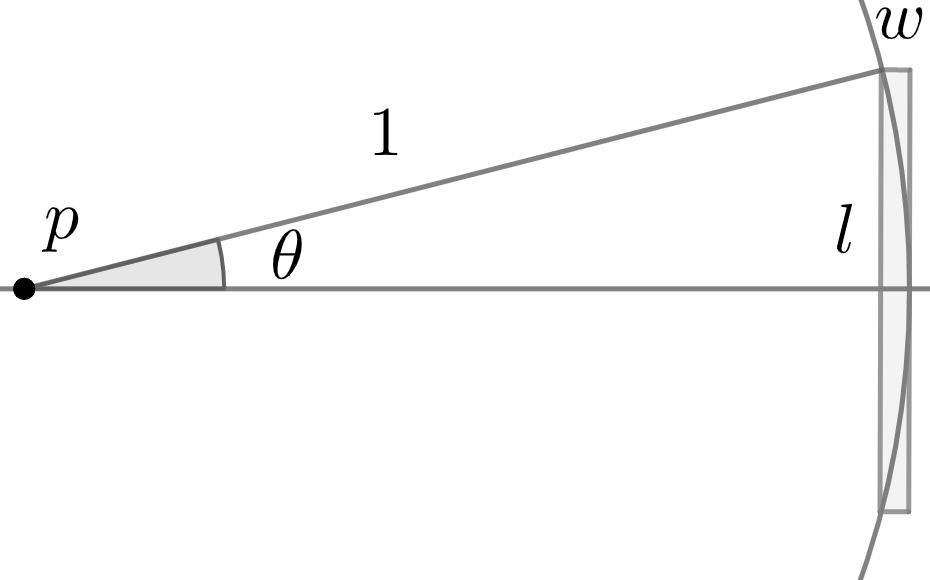}
    \caption{A unit circle centered at $p$ and a circle arc of angle $2\theta$ in a gray rectangle of width $w$ and length $l$.}
    \label{fig:circle_arc_tube}
\end{wrapfigure}

    \noindent \textit{Proof.}
        Let $p$ and $\theta$ be as in Figure \ref{fig:circle_arc_tube} and let the circle arc contained in the gray rectangle be of length $2\sqrt{\delta}$ measured in radians. Thus $\theta=\sqrt{\delta}$. Let $l$ be the length of the rectangle and $w$ be the width. Then $l=2 \sin \theta \approx 2 \sqrt{\delta}$ and $w = 1-\cos \theta \approx \frac{\delta}{2}$. (By small angle approximation). So a shorter arc of length $\sqrt{\delta}$ must be contained in a tube of width $\leq \delta$ tangent to the unit circle.
    \qed

\begin{lem} \label{lem:circlesNminusball}
    Let $(\lines,\pts)$ be an extremal unit circle configuration. Then there exists a refinement $(\lines',\pts')$ such that $\pts'$ is in a ball of radius $N^{-}$ and the set of center points of the unit circles in $\lines'$ is in another ball of radius $N^{-}$.
\end{lem}

\begin{proof}
    Let $\{U_j\}$ be a finite set of disjoint  squares of side length $N^-$ which cover $\pts$. Note that $|\pts|\sim N < \infty$ so we are guaranteed to find a cover of size $\lesssim N$. Let $U^r := \{U_i^r\} \subset\{U_j\}$ such that for all $i$ we have $r \leq |U_i^r \cap \pts|\leq 2r$ and such that $|U^r \cap \pts| \geq \frac{|\pts|}{\log_2N}$. (This is possible to find because there are $\log_2N$ many choices of $r$ and the set of all the $U^r$ cover $\pts$.) Note that a unit circle in $\lines$ goes through at most $N^+$ many squares because the squares are disjoint and have side length $N^-$.
    
    \[ N^{{4 \over 3}} \sim |I(\lines, \pts)| \sim |I(\lines, U^r \cap\pts)| \sim \sum_i |I(U_i^r \cap \lines, U_i^r \cap \pts)| \lesssim \sum_i |U_i^r \cap \lines|^{{2 \over 3}} |U_i^r \cap \pts|^{{2 \over 3}} + |U_i^r \cap \lines| + |U_i^r \cap \pts| \]

    In the last inequality we used Theorem \ref{thm:circle_Sz-trotter} (Szemer\'edi-Trotter for unit circles). Note that in the last sum the contributions from $|U_i^r \cap \pts|$ sum up to $|\pts|\sim N << N^{{4 \over 3}}$ so these terms do not dominate. Furthermore $\sum_i |U_i^r \cap \lines| \lesssim |\lines| N^+ \lesssim N^{1+} << N^{{4 \over 3}}$. This is because the left hand side counts the number of line-square intersections and we recall that each line goes through $\lesssim N^+$ many squares. Thus we can remove the last two terms from the bound. Now we apply H\"older's inequality:

    \[ N^{{4 \over 3}} \lesssim \left( \sum_i |U_i^r \cap \lines| \right)^{{2 \over 3}} \left( \sum_i |U_i^r \cap \pts|^2 \right)^{{1 \over 3}} \lesssim |\lines|^{{2 \over 3}} \left(\frac{|\pts|}{r} r^2 \right)^{{1 \over 3}} \sim |\lines|^{{2 \over 3}} |\pts|^{{1 \over 3}} r^{{1 \over 3}} \sim N r^{{1 \over 3}}\]

    In the second inequality we used $\sum_i |U_i^r \cap \lines| \lesssim |\lines| N^+$ (the $N^+$ was absorbed into the $\sim$). Note this is where we get our scale constraint from: if we chose the squares to have side length much smaller than $N^-$ then each unit circle would intersect $>> N^+$ many squares and this inequality would become correspondingly less sharp.

    Thus $r \gtrsim N$ so there must be some ball of radius $N^-$ which contains $U_i^r$ and therefore contains $|U_i^r \cap \pts| \sim N$ many points from $\pts$.

    We refine the point set to keep only the points inside the ball. We refine the circles to keep only circles with $\gtrsim N^{{1 \over 3}-}$ incidences. (This is a valid refinement by the dual of Lemma \ref{pointsinacellcircle}). The only circles left must have their center point in an annulus of radius 1 and width $N^-$ centered at the ball which contains the points. We can cover this annulus with $N^+$ many finitely overlapping balls of radius $N^-$. One of these balls must contain $\gtrsim N^{1-}$ many center points of unit circles in $\lines$. This yields the desired refinement.
\end{proof}

As before, a key ingredient for our later proofs will be the crossing number inequality which we state here. This is nearly the same result as the standard crossing number inequality from Lemma \ref{crossingnumber} with the key difference that edges are unit circle arcs.

\begin{defn}
    A \textit{crossing} between unit circle arcs is a point in their intersection. 
\end{defn}

Note that a pair of circle arcs can have zero, one, or two crossings. 

\begin{lem}\cite{S} \label{crossingnumbercircle} Let $\pts$ be a set of $n$ points in the plane and $\lines$ be a set of unit circles. Let $G$ be the multigraph where $\pts$ is the vertex set and for every circle $c$ in $\lines$ and pair of adjacent points $(p,p')$ on $c$ we add an edge between $p$ and $p'$. Let $e$ be the number of edges in $G$. Then the number of crossings between edges is $\gtrsim {e^3 \over n^2}$. \end{lem}

\begin{defn}
    We say two circle arcs share a \textbf{simple crossing} if they have exactly one point in their intersection and are not tangent. We say two circle arcs share a \textbf{double crossing} if they have exactly two points in their intersection or if they are tangent. Similarly, two circles have a simple crossing if they each contain an arc bounded by two adjacent points such that these arcs share a simple crossing. Likewise for double crossing.
\end{defn}

Note that when two circle arcs share a double crossing, they are the only pair with a double crossing at that pair of points. 

\begin{lem}\label{simple_crossings}
        Every cell in a cell decomposition from Corollary \ref{structuredcellscircle} has $\gtrsim N^{{4 \over 3}-}$ simple crossings.
    \end{lem}

    \begin{proof}
        There are $\sim N^{{2 \over 3} \pm}$ segments per cell each pair of which contributes $\leq 2$ crossings so each segment participates in $\lesssim N^{{2 \over 3}-}$ many crossings. By Lemma \ref{crossingnumbercircle} there are $\gtrsim N^{{4 \over 3}-}$ crossings so there are $\gtrsim N^{{2 \over 3}-}$ many segments with $\gtrsim N^{{2 \over 3}-}$ crossings each. Assume one of these edges does not contribute $\gtrsim N^{{2 \over 3}-}$ simple crossings. Then it must have $\gtrsim N^{{2 \over 3}-}$ double crossings. By Lemma \ref{lem:circlesNminusball} all edges are circle arcs of length less than $N^{-\epsilon}$ so double crossings correspond to circles being tangent with small error. Quantitatively, given a fixed edge that has double crossings with two other edges, their centers are at distance $\lesssim N^{-\epsilon}$ (small angle approximation). 

        Then also by small angle approximation, the length of the the circle arc spanning the two intersections between circles whose centers are at distance $\lesssim N^{-\epsilon}$ is $\sim 1$ i.e. much larger than $N^{-\epsilon}$. Thus no edge can span both the crossings of a pair of circles with close centers. So if two edges form double crossings with a third edge, then the first two must form a simple crossing. Furthermore this simple crossing must be contained in the same cell since all edges are entirely contained in a cell.

        So if one of these crossing rich edges $e$ in a cell does not contribute $\gtrsim N^{{2 \over 3}-}$ simple crossings, then every pair of edges that form double crossings with $e$ must simple cross each other. But $e$ is crossing rich so it must have $\gtrsim N^{{2 \over 3}-}$ double crossings which implies $\gtrsim N^{{4 \over 3}-}$ simple crossings in the cell. Thus every cell has $\gtrsim N^{{4 \over 3}-}$ simple crossings.
    \end{proof}

\subsection{Structuring Circles, Organizing Circles and the Bush Construction}\label{subsection:structuring_organizing_bush_circle}

\begin{wrapfigure}[18]{R}{0.3\textwidth}
    \centering
    \includegraphics[width=0.3\textwidth]{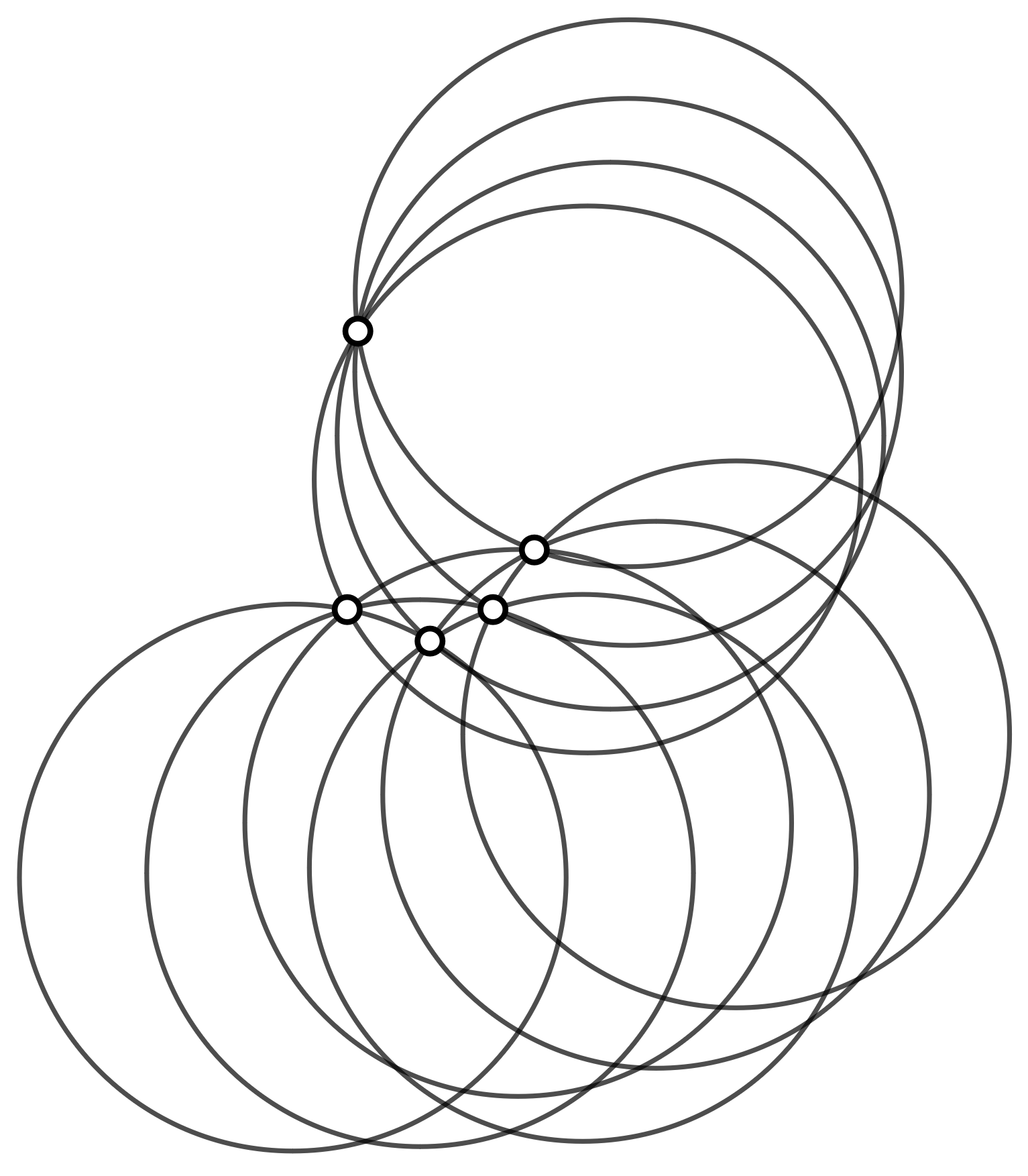}
    \caption{Five structuring points structure the set of black unit circles. }
    \label{fig:structured_circles}
\end{wrapfigure}

We make a definition of a structured set of circles.

\begin{defn}  We say that a set $\lines_1$ of at least $N^{{2 \over 3}-}$ unit circles is \textbf{structured} if there is
a set $\pts_1$ of at most $N^{{1 \over 3}+}$ points so that each circle of $\lines_1$ is incident to at
least two points of $\pts_1$. We call this set of points \textbf{structuring}. See Figure \ref{fig:structured_circles}. \end{defn}

Note that since $\lesssim N^{\frac{2}{3}+}$ circles go through at least two among $\lesssim N^{\frac{1}{3}+}$ points, the structuring points essentially define the structured circles. Now we're ready to state our structuring theorem.

\begin{figure}[!h]
    \centering
    \includegraphics[width=\textwidth]{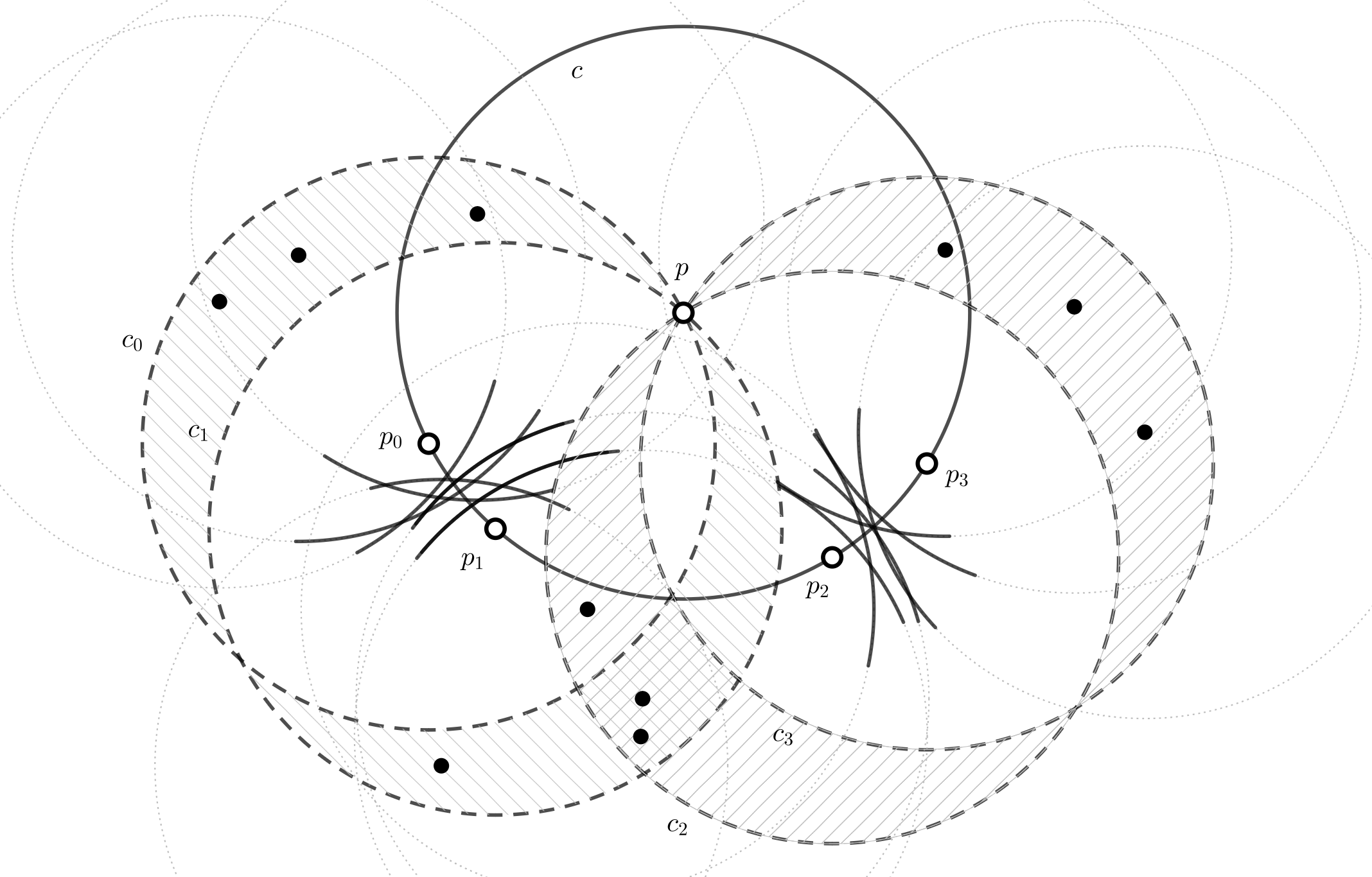}
    \caption{An organizing circle $c$ has two circle arcs bounded by $p_0,p_1$ and $p_2,p_3$ respectively, where both arcs are crossed by a structured set of circles (dotted light gray with crossing arcs in solid black). The dual formulation is shown: the centers of the structured circles form a set of structured points (black) contained in the interior of the sectors defined below (hashed gray region bounded by dotted black circles which intersect at the organizing point $p$: center of $c$).}
    \label{fig:organizing_circle}
\end{figure}

\begin{thm} \label{verynicespacingcircle}  Let $(\lines,\pts, \cal{J})$ be an extremal partial unit circle configuration. Then there exists a point weighted cell decomposition with $\gtrsim N^{{2 \over 3}-}$ cells and a refinement of the configuration $(\lines^{\prime},\pts,\cal{J}')$ so that for each circle $c \in \lines^{\prime}$ there are points $p_0, \dots
, p_M$ of $\pts$ with $(c,p_j) \in \cal{J}'$ and the $p_j$'s in order of their position on $c$ and with $M \gtrsim N^{{1 \over 3}-}$ so that
for each consecutive pair of points $p_{2k},p_{2k+1}$, the circle arc bounded by $p_{2k},p_{2k+1}$ is entirely contained in a cell, and there is a structured set of circles $\lines_k$ so that each $c^{\prime}$ in $\lines_k$ intersects $c$ in the open circle arc bounded by the points $p_{2k}$ and $p_{2k+1}$ exactly once. We say the circles in $\lines^{\prime}$ \textbf{organize} $\pts$. See Figure \ref{fig:organizing_circle}. \end{thm}

Note that we restrict our attention to crossings where the circle arcs share only one point in common instead of two.

\begin{proof}   First consider an extremal partial unit circle configuration $(\lines,\pts, \cal{J})$. We consider the point weighted cell decomposition $C_1, \dots , C_{N^{{2 \over 3}-}}$ guaranteed by Theorem \ref{nicerefinementcircle} and Corollary \ref{structuredcellscircle}. We obtain the refined unit circle incidence multigraph $G$ where the vertex set is $\pts$ and for every pair of vertices, we add an edge if the corresponding points are adjacent on a unit circle, and the unit circle arc that they bound has length $\lesssim N^{-\epsilon}$ and is entirely contained inside a cell and such that each cell has $\lesssim N^{{2 \over 3}+}$ edges. By Corollary \ref{structuredcellscircle} we know each cell has $\gtrsim N^{{2 \over 3}-}$ edges of $G$ and $\lesssim N^{{1 \over 3}+}$ points.

We know from Lemma \ref{simple_crossings} that there are $\gtrsim N^{{4 \over 3}-}$ simple crossings per cell. Each cell has $\lesssim N^{{2 \over 3}}$ circle arcs (because it has $\lesssim N^{{1 \over 3}}$ points and each pair of points is connected by at most 2 circle arcs) so each cell has $\gtrsim N^{{2 \over 3}-}$ circle arcs that each have $\gtrsim N^{{2 \over 3}-}$ simple crossings. We call these \textbf{sector arcs} because their $\gtrsim N^{{2 \over 3}-}$ simple crossings come from a set of $\gtrsim N^{{2 \over 3}-}$ structured circles which are structured by the set of points in the cell.

There are $\gtrsim N^{{2 \over 3}-}$ cells so $\gtrsim N^{{4 \over 3}-}$ sector arcs total. There are $\lesssim N$ circles in the configuration each of which (under our refinements) has no more than one circle arc per cell so has $\lesssim N^{{1 \over 3}}$ sector arcs. Thus $\gtrsim N^{1-}$ circles must have $\gtrsim N^{{1 \over 3}-}$ sector arcs.
\end{proof}

We'd now like to take advantage of our result. Theorem \ref{verynicespacingcircle} is a method of associating to each extremal unit circle configuration
a refinement which is rather nicely parametrized. 
It gives us many circles $c$ which are incident to particular sets of points $p_1, \dots ,p_M$ with $M \gtrsim N^{{1 \over 3}-}$
so that we have $\gtrsim N^{{2 \over 3}-}$ circles intersecting $c$ between adjacent points which are structured. We now apply point--circle duality which maps a point to the unit circle centered at it and vise-versa, obtaining a dual form of Theorem \ref{verynicespacingcircle}.

\begin{defn}
    Consider a point $p$ with a bush $c_1, \dots, c_M$ of circles going through it ordered according to the direction of their tangent at $p$. Let $d_1, \dots, d_M$ be the interiors of the circles $c_i$. This bush of circles defines $M$ \textbf{sectors} where we define the $i^{th}$ sector $s_i = (d_i \cap d_{i+1}^c) \cup (d_i^c \cap d_{i+1})$ where $d^c$ is the complement of the disk $d$.
\end{defn}

We analyze what it means for a point to be in the interior of a circle in terms of circle crossings.

\begin{wrapfigure}[23]{r}{0.65\textwidth}
    \centering
    \includegraphics[width=0.65\textwidth]{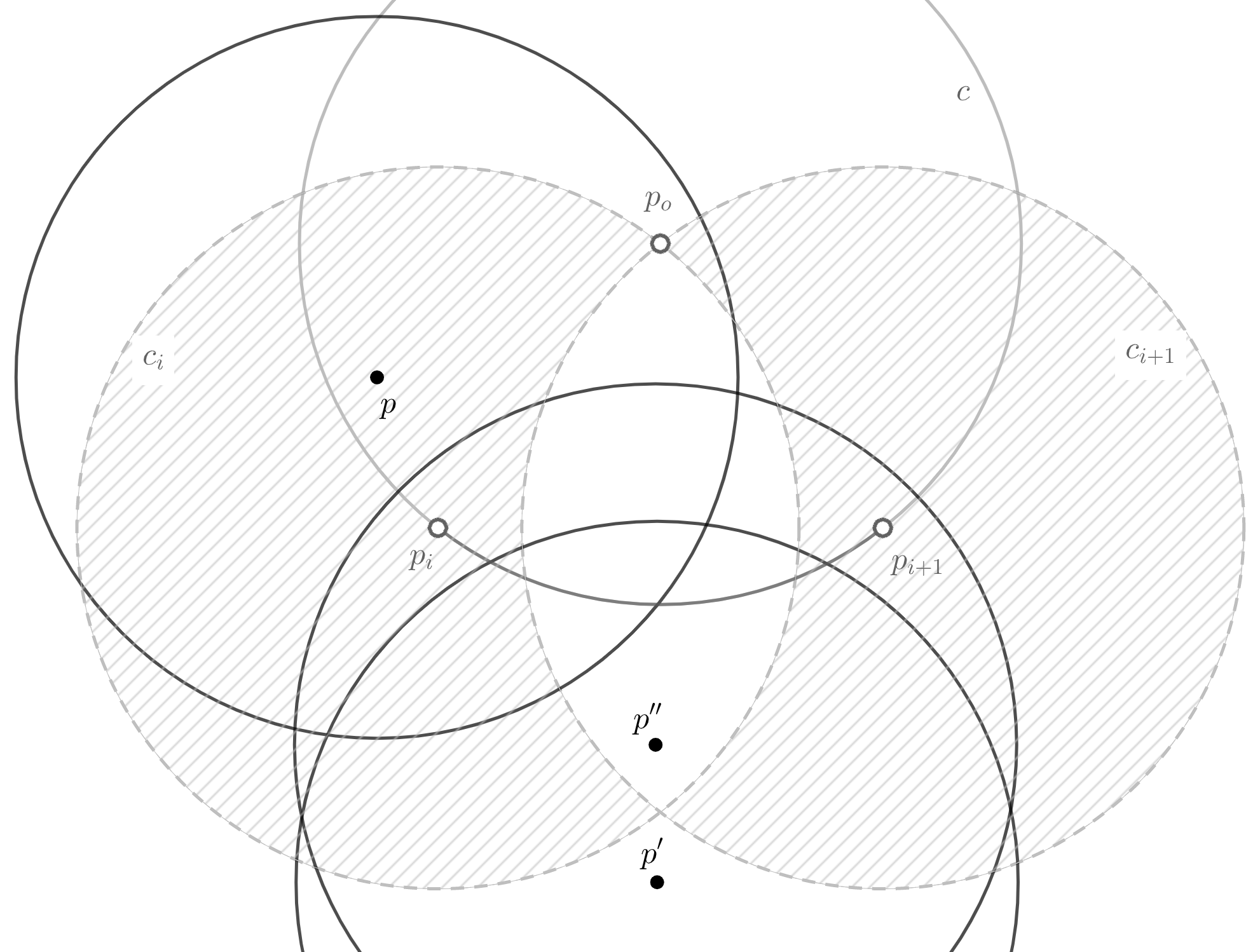}
    \caption{An organizing circle $c$ and its sector $s_i$ (hatched area) are shown in gray. $p$ is in $s_i$ but $p'$ and $p''$ are not. The circle centered at $p$ crosses the circle arc of $c$ between $p_i,p_{i+1}$ exactly once. The circle centered at $p'$ crosses it twice, and the circle centered at $p''$ crosses it zero times.}
    \label{fig:circle_sectors}
\end{wrapfigure}

\begin{lem}\label{circleduality}
    A point $p$ is in the sector $s_i$ of an organizing circle $c$ if and only if the unit circle centered at $p$ crosses the arc of $c$ between $p_i$ and $p_{i+1}$ exactly once. See Figure \ref{fig:circle_sectors}.
\end{lem}

\noindent \textit{Proof.} Given two circles $c$ (shown in gray in Figure \ref{fig:circle_sectors}) and $c'$ (three possible examples are shown in black in Figure \ref{fig:circle_sectors}) and points $p_i$ and $p_{i+1}$ on $c$, $c'$ crosses the circle arc of $c$ contained between $p_i$ and $p_{i+1}$ exactly once if and only if exactly one of $p_i$ or $p_{i+1}$ is in the interior of $c'$. This happens if and only if $p$, the center of $c'$, is in exactly one of the interiors of $c_i$ or $c_{i+1}$, the circles centered at $p_i$ and $p_{i+1}$ respectively. Note that being in the symmetric difference of the disk bounded by $c_i$ and of the disk bounded by $c_{i+1}$ is equivalent to being in the sector $s_i$. Thus $c'$ crosses the circle arc of $c$ contained between $p_i$ and $p_{i+1}$ exactly once if and only if the point $p$ (center of $c'$) is in the sector $s_i$ of the organizing point $p_o$, where $p_o$ is the center of $c$.
\qed

\bigskip

This motivates the our need to distinguish between \textit{simple crossings} and \textit{double crossings}. Only simple crossings tell us information about points/ circles belonging to sectors. Double crossings instead tell us that two circles are almost tangent but this only happens for a very small fraction of circle crossings (see Lemma \ref{simple_crossings}). So double crossings with an organizing circle can only partition a small fraction of circles from the example which is not useful for our purpose. We obtain the dual of Theorem \ref{verynicespacingcircle}:

\begin{thm} \label{thm:dualverynicespacingcircle}  Let $(\lines,\pts, \cal{J})$ be an extremal partial unit circle configuration. Then there is a refinement $(\lines,\pts^{\prime}, \cal{J}')$ so that for each point $p \in \pts^{\prime}$ there are circles $c_0, \dots 
,c_M$ of $\lines$ which are incident to $p$ in $\cal{J}'$ in order of their direction and with $M \gtrsim N^{{1 \over 3}-}$ so that each sector $s_j$ bounded by consecutive pairs of circles $c_{2j},c_{2j+1}$ contains a structured set of $\gtrsim N^{\frac{2}{3}-}$ points $\pts_j$. We say the points in $\pts^{\prime}$ \textbf{organize} $\lines$. \end{thm}

\begin{proof}
    Let $D(\cdot)$ be the point-unit circle dual operation that takes in a point in $\RR^2$ and outputs the unit circle centered at that point and takes in a unit circle in $\RR^2$ and outputs the center point. Let $D(\cal{J})$ be the set of incidences between points $D(c)$ and unit circles $D(p)$ such that $(p,c)\in \cal{J}$. By Lemma \ref{circleduality}, the statement of Theorem \ref{thm:dualverynicespacingcircle} applied to the partial configuration $(\lines,\pts, \cal{J})$ is equivalent to the statement of Theorem \ref{verynicespacingcircle} applied to the partial configuration $(D(\pts), D(\lines), D(\cal{J}))$. 
    
\end{proof}

For any point $p \in \pts^{\prime}$ with circles
$c_1, \dots c_M$ incident to it, there are $\gtrsim N^{{2 \over 3}-}$ points in each of the $\gtrsim N^{{1 \over 3}-}$ sectors for a total of $N^{1-}$ points. We take this set of points as a refinement $\pts^{\prime}$ of our original set of points $\pts$. What is particularly pleasant about this structure is that each of the $\gtrsim N^{{2 \over 3}-}$ points
of $\pts^{\prime}$ in a sector $s_j$ lie on at least two of the $\gtrsim N^{{1 \over 3}-}$ structuring circles.

Structuring circles seem very odd precisely because all of the points on them lie in a particular sector $s_j$. But this is
not as odd as it seems. Every cell has $\approx N^{{1 \over 3} \pm}$ points in it and $\approx N^{{2 \over 3} \pm}$ circles going through it so by duality the set of $\approx N^{{1 \over 3} \pm}$ structuring circles define the $N^{{2 \over 3} \pm}$ points in a sector.

We're going to show that for any choice of organizing point $p \in \pts^{\prime}$ a typical circle $c$ will have incidences in most sectors $s_j$. To do this we first need to introduce a refinement of the configuration endowed with a cell decomposition whose boundary circles include the bush through $p$. 

\begin{thm}[bush construction] \label{bushconstructioncircle}
For any extremal unit circle configuration $(\lines, \pts)$ there exists a subset $\pts^{\prime}$ of $\gtrsim N^{1-}$ points in $\pts$ which are organizing with $\sim N^{{1 \over 3} \pm}$ sectors and a refined configuration $(\lines^{\prime}, \pts^{\prime})$ such that the $\gtrsim N^{1-}$ circles in $\lines^{\prime}$ organize $\pts^{\prime}$. Also for any $p\in \pts^{\prime}$ the refinement $(\lines^{\prime}, \pts_p)$ where $\pts_p$ are the points in $\pts$ organized by $p$, has a refinement
$(\lines^{\prime}, \pts_p^{\prime})$ which is an extremal
configuration with a point-weighted cell decomposition where each cell is contained in a sector. Moreover, any circle $c \in \lines^{\prime}$ which has exactly $N^{{2 \over 3}+\alpha}$
simple crossings with circles of $\lines^{\prime}$ within the sector $s$ with $\alpha>k\epsilon$ for $k$ sufficiently large
will not enter more than $N^{\alpha+}$ cells in $s$.
 \end{thm}

 Note as before we split our circles into circle arcs between adjacent points so we say two circles have a simple crossing if they each have a circle arc between a pair of adjacent points such that the intersection of these circle arcs is exactly one point and the arcs are not tangent.

\begin{proof}

\textit{Refinement properties:} We apply Theorem \ref{thm:dualverynicespacingcircle} and obtain the refinement $(\lines, \pts^{\prime})$ of organizing points. Then we apply Theorem \ref{verynicespacingcircle} to $(\lines, \pts^{\prime})$ obtaining the refinement $(\lines^{\prime}, \pts^{\prime})$ where the circles in $\lines^{\prime}$ organize $\pts^{\prime}$. Now we have shown the first claim of the theorem.

\bigskip


\textit{Cell decomposition:} Let $p \in \pts^{\prime}$ and $\pts_p$ be the set of points in $\pts$ organized by $p$. Note there are $\gtrsim N^{1-}$ organized points each $\gtrsim N^{{1 \over 3}-}$ rich so $(\lines^{\prime}, \pts_p)$ is an extremal configuration which we work in for this paragraph. We label the bush of $M \sim N^{{1\over 3}\pm}$ circles intersecting $p$ as $c_0,\dots, c_M$.
By Theorem \ref{nicerefinementcircle} we can find a set of unit circles $\lines_{N^{{1 \over 3}}} \subset \lines'$ which yields a point weighted cell decomposition for some refinement of the point set $\pts_p'$. Our ``cell decomposition" will be made from
the circles $c_0,\dots,c_M$ together with the circles in $\lines_{N^{{1 \over 3}}}$. 

But we will not think of this collection of circles as giving a cell decomposition in the conventional sense. Recall that the set of points $E_c$ in $\RR^2$ which are the centers of
circles that simple cross $c_0$ are double-covered by the sectors  $s_i = (d_i \cap d_{i+1}^c) \cup (d_i^c \cap d_{i+1})$ defined by the adjacent circles $c_i$ and $c_{i+1}$. For
each sector $s_i$ we let the circles $c^{\prime}_1, \dots,c^{\prime}_K$ subdivide
$s_i$ into cells. We have now obtained a collection of cells which double cover $E$. This
will essentially be all right for us. We will sometimes over count incidences, but we
will at most double count them.

The cell decomposition bounded by
the circles $c_0,\dots,c_M$ together with the circles in $\lines_{N^{{1 \over 3}}}$ is still point weighted because the decomposition with just $\lines_{N^{{1 \over 3}}}$ was point weighted and adding in more boundary circles can only decrease the number of points per cell. Also we added in $\lesssim N^{{1 \over 3}+}$ boundary circles so generic circles still enter at most $\lesssim N^{{1 \over 3}+}$ cells. Thus this decomposition is point weighted and each cell is contained in one sector.

\bigskip

\textit{Crossings} We use the second part of the $\lines_{N^{1 \over 3}}$ partitions circle-boundary arc intersections property in Lemma \ref{lem:boundary_circles}. This property says that for any $c \in \lines$ and any list of $\log(N)N^{2 \over 3}$ unit circles which cross $c$ consecutively, at most $10 \log(N)$ of them are in $\lines_{N^{1 \over 3}}$. So inside a sector $s$, $c$ takes $\gtrsim N^{{2 \over 3}-}$ crossings per $10 \log(N)$ consecutive cells that $c$ intersects. Thus if $c$ has $\sim N^{{2 \over 3}+\alpha}$ crossings in $s$, then $c$ enters $\lesssim N^{\alpha+}$ cells in $s$.

 \end{proof}

When trying to show that a typical circle in a configuration $(\lines,\pts)$ will have incidences in most of the sectors of an organizing point, the enemy case is circles which take too many points in a given sector. By Theorem \ref{boundedsidescircle} we know circles take $\lesssim N^+$ incidences per cell so the enemy case is circles that have incidences from $N^{\alpha+}$ cells in the sector where $\alpha > k\epsilon$. If we show that circles with $N^{{2 \over 3}+\alpha}$ simple crossings in a sector, which we will call \textit{$\alpha-$fast circles}, do not contribute significantly to the total number of incidences in that sector, then most incidences come from circles with $\lesssim N^{{2 \over 3}+}$ simple crossings in that sector and by Theorem \ref{bushconstructioncircle} we know these circles enter $\lesssim N^{+}$ cells in that sector so we win.

\begin{thm} \label{circlesectorincidences}
There exist $\gtrsim N^{1-}$ organizing points $p$ in $(\lines, \pts)$ with $\sim N^{{1 \over 3}\pm}$ sectors such that for every sector $\gtrsim N^{1-}$ of its incidences come from circles taking $\lesssim N^{-}$ points in that sector.
\end{thm}

\noindent\textit{Proof.} Let $\epsilon$ be such that $I(\lines,\pts)\gtrsim N^{{4 \over 3}-\epsilon}$ and $|\pts|, |\lines| \gtrsim N^{1 -\epsilon}$ and such that the $\pm$ in the statements of Theorems \ref{thm:dualverynicespacingcircle} and \ref{bushconstructioncircle} have magnitude smaller than $\epsilon$. We  refine the point set to keep only points which are $\gtrsim N^{{1 \over 3}-\epsilon-}$ rich. ($|\pts| \lesssim N$ so we removed at most $N^{{4 \over 3}-\epsilon-}$ incidences.) We further refine the point set to keep only points that are $\lesssim N^{{1 \over 3}+\epsilon}$ rich. Note that Theorem \ref{thm:circle_Sz-trotter} (Szemer\'edi-Trotter theorem for circles) implies this last refinement removed $\lesssim N^{{4 \over 3}-2\epsilon}$ incidences so this is a valid refinement. Next choose an organizing point $p$ and use the bush construction $(\lines', \pts_p')$ from Theorem \ref{bushconstructioncircle}. Note the bush at $p$ has $\gtrsim N^{{1 \over 3} -\epsilon-}$ many circles.

\bigskip

We first remove sectors where the pair of lines tangent to the sector's boundary circles which go through the organizing point make an angle $ \gtrsim N^{-{1 \over 3} +2\epsilon +}$. These angles partition the unit circle so we removed $ \lesssim N^{{1 \over 3} -2\epsilon -}$ many sectors and kept only $\textit{narrow sectors}$.

\bigskip

Next we remove sectors that have $\gtrsim N^{{5 \over 3}+2\epsilon}$ simple crossings. There are at most $N^{{1 \over 3}-2\epsilon}$ such sectors because $\lesssim N^{2}$ simple crossings total because each pair of circle contributes at most two simple crossings. Similarly we remove sectors that have $\lesssim N^{{5 \over 3}-2\epsilon}$ simple crossings. So we kept $\gtrsim N^{{1 \over 3}-\epsilon}$ sectors with $\sim N^{{5 \over 3} \pm}$ crossings. Also note from the definition of organizing point in Theorem \ref{thm:dualverynicespacingcircle} that each sector has $\gtrsim N^{{2 \over 3}-}$ points which from our first refinement are all $\gtrsim N^{{1 \over 3}-}$ rich so each sector contributes $\gtrsim N^{1-\epsilon}$ incidences.

\bigskip

Let this valid refinement be $(\lines^{\prime}, \pts^{\prime})$. By Corollary \ref{structuredcellscircle} we may chose a subset $J(\lines^{\prime}, \pts^{\prime}) \subset I(\lines^{\prime}, \pts^{\prime})$ such that $|J(\lines^{\prime}, \pts^{\prime})| \gtrsim N^{{4 \over 3}-}$ and every circle has $\lesssim N^{+}$ incidences from $J(\lines^{\prime}, \pts^{\prime})$ per cell.

\begin{defn}[Fast circles]    \label{fastlinescircle}  We say that a circle is $\alpha-$fast for a sector $s$ if it has $\sim N^{{2 \over 3}+\alpha\pm}$ simple crossings with circles in $\lines$.
\end{defn}

This is the enemy case.

\bigskip

\textbf{Claim:} \textit{An $\alpha-$fast circle goes through $\sim N^{\alpha \pm}$ cells in $s$.}

The boundary circles in $s$ came from a cell decomposition that is guaranteed by Lemma \ref{lem:boundary_circles} to satisfy good partitioning properties, including the first partitions intersections between unit circles and boundary arcs property. Thus an $\alpha-$fast circle (which crosses $N^{{2 \over 3}+\alpha \pm}$ circles of $\lines$ in $s$) must cross at least $N^{\alpha-}$ of the cell boundary circles in the sector and
therefore must enter at least $N^{\alpha-}$ cells. Similarly by Theorem \ref{bushconstructioncircle} each $\alpha-$fast circle enters no more than $N^{\alpha+}$ cells in $s$. 

\bigskip

Let $c_1$ be the first boundary circle of the sector $s$ and let $c_2$ be the second boundary circle. $s$ is the union of two connected components $s_1$ and $s_2$. The sector $s$ has $\gtrsim N^{1-\epsilon}$ incidences so WLOG let $s_1$ be a connected component that has $\gtrsim N^{1-\epsilon}$ incidences. We will show that fast circles do not contribute significantly to the number of incidences in $s_1$. 

\bigskip

\textbf{Claim} There exists a structuring circle $c$ for the sector $s$ which takes $\gtrsim N^{{1 \over 3}-4\epsilon}$ incidences in $s_1$. 

Proof: If each structuring circle from $s$ has $\lesssim N^{{1 \over 3}-4\epsilon}$ points in $s_1$ then there are $\lesssim N^{{2 \over 3}-3\epsilon+}$ points in $s_1$ since there are $\lesssim N^{{1 \over 3}+\epsilon}$ many structuring circles. Each point is $\lesssim N^{{1 \over 3} +\epsilon}$ rich (see refinement at beginning of the proof) so this implies $s_1$ has $\lesssim N^{1-2\epsilon+}$ incidences.

\bigskip

Also this structuring circle $c$ has $\lesssim N^{+}$ incidences from $\cal{J}(\lines', \pts')$ per cell so $c$ goes through $\gtrsim N^{{1 \over 3}-\delta}$ cells in $s_1$ (where $\delta$ can be computed from the exponents in the statements of Theorems \ref{thm:dualverynicespacingcircle} and \ref{boundedsidescircle} and we require $\delta > 4\epsilon$, independent of the choice of $s$). We do a refinement of $s_1$ where we only keep points in cells that $c$ has an incidence with. 

\bigskip

We refine the boundary circles in $s$ by removing all the $\alpha-$fast ones for $\alpha > 4\epsilon$. The number of fast circles is $\lesssim N^{1-4\epsilon}$ since there are $\lesssim N^{1+}$ simple crossings in $s$ and each fast circle takes $\gtrsim N^{4\epsilon}$ many. So the event that at least half the boundary circles are fast occurs with probability $\lesssim N^{-4\epsilon N^{{1 \over 3}-}}$. We may exclude this unlikely event when selecting our boundary circles from Lemma \ref{lem:boundary_circles}. Thus $c$ still intersects $\gtrsim N^{{1 \over 3}-4 \epsilon}$ boundary circles after our refinement. The remaining boundary circles partition $c$ into $\gtrsim N^{{1 \over 3}-4\epsilon}$ many disjoint circle arcs. Each circle arc from $c$ contained in a cell must be adjacent to an intersection of $c$ with a boundary arc of the cell. Each boundary arc intersects $c$ at most twice and by Theorem \ref{nicerefinementcircle} each cell has $\lesssim N^+$ boundary arcs. So each cell contains $\lesssim N^+$ many circle arcs from $c$. Thus we kept $\gtrsim N^{{1 \over 3}-4\epsilon}$ cells which $c$ goes through so this is still a valid cell decomposition.

\bigskip

Now we show that circles which cross $s_1$ must intersect both $c_1$ and $c_2$ in $s_1$. This will allow us to order them by where they enter and exit $s_1$. Recall from the beginning of the proof that we kept only narrow sectors. If a circle $c'$ enters $s_1$ but does not intersect both $c_1$ and $c_2$ then $c'$ must be $\lesssim N^{-{1 \over 3}+}$ away from being tangent to $c_1$. This implies that the centers of $c_1$ and $c'$ are distance $>\frac{1}{2}$. We apply Lemma \ref{lem:circlesNminusball} to refine the circles such that all remaining circles are distance $\lesssim N^-$ from each other so this almost tangency situation cannot occur. Thus any circle which crosses $s_1$ must intersect both $c_1$ and $c_2$ in $s_1$.

\bigskip

Let $i_{1,j}$ be the set of $\alpha-$fast circle arcs in $s_1$ which intersect $c_1$ in the interval bounded by its $j^{th}$ and $(j+1)^{th}$ intersection points with the remaining boundary circles. Similary define $i_{2,j}$ in terms of the intersection order with $c_2$. Note that if $c' \in i_{1,j} \cap i_{2,k}$ then $c'$ intersects at least $|j-k|$ many cells because $c'$ must intersect at least $j-k$ many boundary circles which intersect $c_1$ above $i_{1,j}$ and intersect $c_2$ below $i_{2,k}$ (with the orders reversed if $k>j$). Note a circle intersects $s_1$ at most twice so a circle belongs to at most two $i_{1,j}$. Now define \[I_{1,j}=\cup_{k \in [ jN^{\alpha},  (j+1)N^{\alpha}]} i_{1,k} \text{ for } j \in [1, N^{{1 \over 3}-\alpha}]. \]

Define $I_{2,j}$ similarly. Since $\alpha-$fast circles go through $\sim N^{\alpha \pm}$ many cells, they cross $\lesssim N^{\alpha+}$ many boundary circles so $I_{1,j} \cap I_{2,k} = \emptyset$ if $|j-k| \geq N^+$.

\bigskip

Now we prove the set of circles in $I_{1,j}$ intersect $\lesssim N^{\alpha}$ many cells in $s_1$. This follows from the claim that any boundary circle which intersects some circle from $I_{1,j}$ must intersect $c_1$ in the interval spanned by $I_{1,j}$, or intersect $c_2$ in the interval spanned $I_{2,j \pm N^+}$. This is because the remaining boundary circles are slow so for $\alpha > \delta$ there are no boundary circles that intersect $c_1$ above $I_{1,j}$ and $c_2$ below $I_{2,j\pm N^+}$ or vice-versa. So for a fixed $j$, there are $\lesssim N^{\alpha+}$ many boundary circles which contribute to cells that circles in $I_{1,j}$ could intersect. The structuring circle $c$ intersects each of these boundary circles at most twice and all the remaining cells are skewered along $c$. So $c$ goes through at most two cells per boundary circle and the circles in $I_{1,j}$ intersect $\lesssim N^{\alpha+}$ many cells, so take incidences with $\lesssim N^{{1 \over 3}+\alpha+}$ many points.

\bigskip

Now we apply Szemer\'edi-Trotter (Theorem \ref{thm:circle_Sz-trotter}) to count the number of incidences from $\alpha-$fast circles by summing over each $I_{1,j}$. Let $\lines_\alpha$ be the set of $\alpha-$fast circles.

\[|I(\pts, \lines_\alpha)| \lesssim \sum_j |I_{1,j}|^{{2 \over 3}} (N^{{1\over 3}+\alpha+})^{2 \over 3} \lesssim |\lines_\alpha|^{2 \over 3} (\sum_j (N^{{1\over 3}+\alpha+})^2)^{1 \over 3} \lesssim N^{(1-\alpha +){2 \over 3}} (N^{{1 \over 3} - \alpha} N^{{2 \over 3}+ 2\alpha+})^{1 \over 3} \lesssim  N^{1-\alpha/3+}\]

The first inequality follows from the Szemer\'edi-Trotter applied to each interval, the second step follows from H\"older's inequality and $\sum_j |I_{1,j}| = |\lines_\alpha| $. To deduce the third inequality we note there are $\lesssim N^{{5 \over 3}+}$ simple crossings in $s$ and each $\alpha-$fast circle takes $\sim N^{{2 \over 3}+\alpha}$ simple crossings so $|\lines_\alpha| \lesssim N^{1-\alpha+}$ and we recall $j \lesssim N^{{1 \over 3} - \alpha}$.

\bigskip

The number of incidences in $s_1$ is $\gtrsim N^{1-\epsilon}$. We showed that $|I(\pts, \lines_\alpha)| \lesssim N^{1-\alpha/3}$ so the total number of incidences from fast circles is $\lesssim N^{1-\delta}<<N^{1-\epsilon}$. Thus the slow circles contribute the majority of the incidences

\qed

\begin{wrapfigure}[28]{R}{0.5\textwidth}
    \centering
    \includegraphics[width=0.5\textwidth]{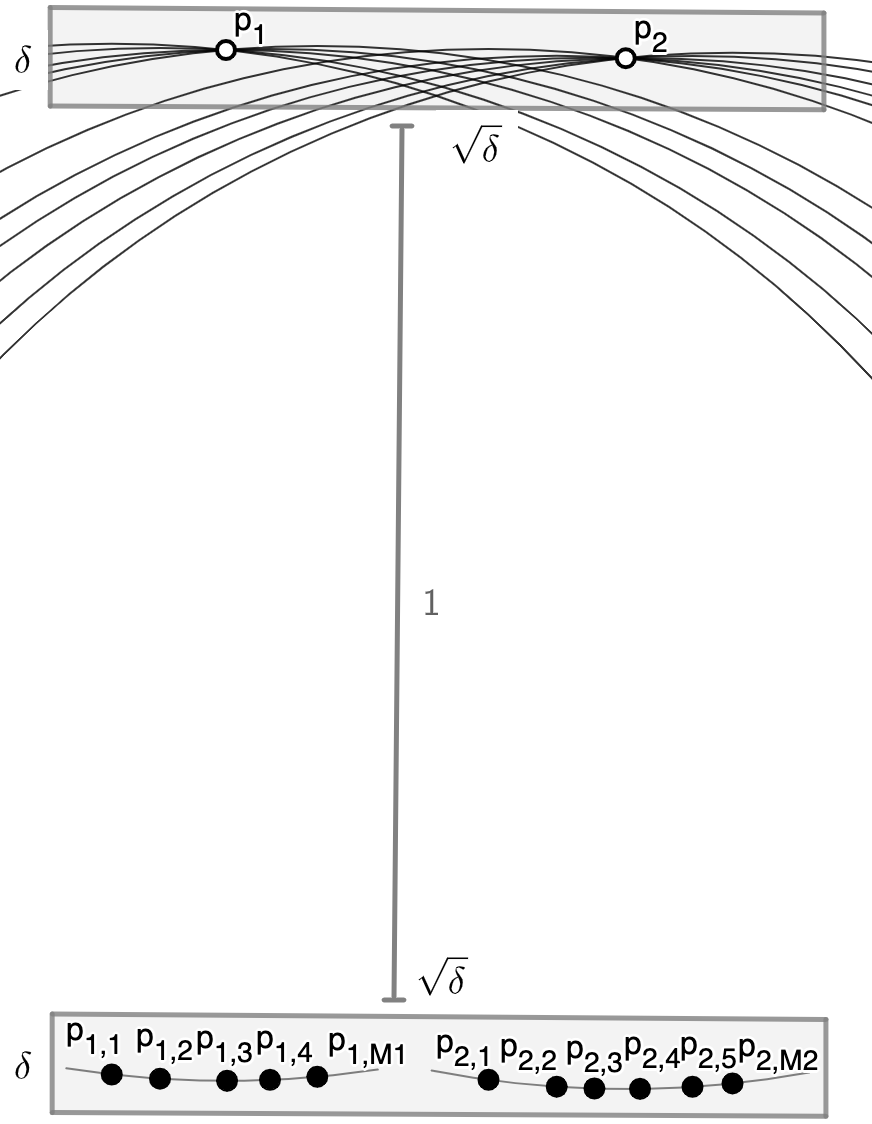}
    \caption{A two bush cell decomposition organized by $p_1$ and $p_2$ with bushes respectively centered at $p_{1,1},...,p_{1,M_1}$ and $p_{2,1},...,p_{2,M_2}$.}
    \label{fig:twobushcircles}
\end{wrapfigure}

\begin{figure}[!h]
    \centering
    \includegraphics[width=1\textwidth]{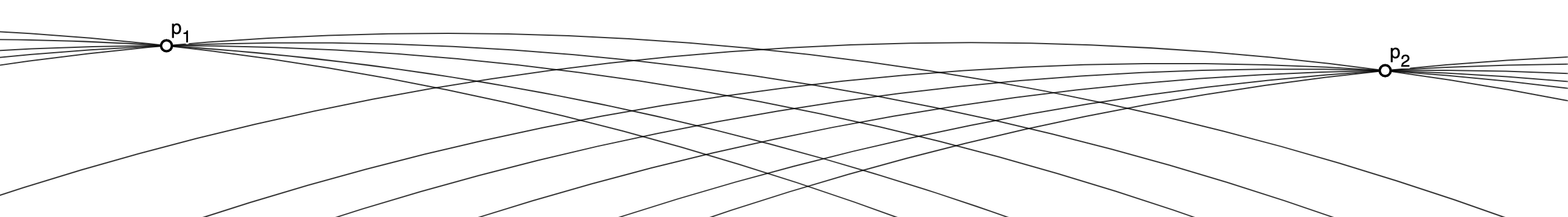}
    \caption{A zoomed in view of the gray box at the top of Figure \ref{fig:twobushcircles} shows the boundary circles decomposing the space into a grid of $N^{2 \over 3}$ many cells. This is the circle analogue of Figure \ref{fig:2bushcelldecomp}.}
    \label{fig:zoomedtwobushcircles}
\end{figure}

\begin{thm} \label{doublebushmixingcircle}
Let $(\lines,\pts)$ be an extremal unit circle configuration.  Then there are $\gtrsim N^{2-}$ pairs of organizing points $p_1$ and $p_2$, with bushes $c_{1,1},\dots, c_{M_1,1}$ incident to $p_1$ and
$c_{1,2}, \dots, c_{M_2,2}$ incident to $p_2$ with $M_1,M_2 \sim N^{{1 \over 3}\pm}$
and a refinement $\pts^{\prime} \subset \pts$ so that the two bushes break 
$\pts^{\prime}$ into $M_1 M_2$ cells which are point weighted. Moreover each sector $s$
of the bush at $p_1$ with at least $N^{{2 \over 3}-}$ points of $\pts^{\prime}$ has
at least $N^{1-}$ circles of $\lines$ incident to at least two points in some cell of $s$. See Figures \ref{fig:twobushcircles} and \ref{fig:zoomedtwobushcircles}.
\end{thm}

\noindent \textit{Proof.} Apply Theorem \ref{verynicespacingcircle}, we find a point $p_1$ with bush
$c_{1,1},\dots, c_{M,1}$ and structuring circles holding in total $N^{1-}$ points. We
call this set of points $\pts_1$. Then $(\lines, \pts_1)$ is an extremal configuration.
We apply Theorem \ref{circlesectorincidences} to find a refinement of the set of incidences
$J(\lines,\pts_1)$ with $|J(\lines,\pts_1)| \gtrsim N^{{4 \over 3}-}$ so that each circle of $\lines$
takes only $N^{+}$ incidences of $J(\lines,\pts_1)$ in each sector of the bush at $p_1$.

We restrict to those points in $\pts_1$ which are at least $N^{{1 \over 3}-}$ rich in incidences of $J(\lines,\pts_1)$.  We refer to that set as $\pts_2$. The set $(\lines,\pts_2)$
is an extremal configuration.  We refine the set of circles to $\lines_1$ which take $N^+$
incidences in $N^{{1 \over 3}-}$ sectors of the bush at $p_1$. We let $\pts_3$ be
the set of points of $\pts_2$ that are $N^{{1 \over 3}-}$ rich with respect to circles of $\lines_1$ incident to only $N^+$ other points in the same sector. The pair $(\lines_1,\pts_3)$ is an extremal configuration. Pick an organizing point $p_2$ with bush 
$c_{1,2} \dots c_{M_2,2}$ and structuring circles for each sector of the bush from $\lines_1$.
Keep only the points on the structuring circles of the second bush which occur in groups of at most $N^{+}$ many per sector of the first bush. Call this set of points $\pts^{\prime}$. There are $\lesssim N^{{1 \over 3}+}$ structuring circles in each sector of the second bush so the cell decomposition given by the two bushes is point weighted.

Since each sector of the first bush has $N^{1-}$ incidences coming from circles making at
most $N^{+}$ incidences in the sector, there must be $N^{1-}$ circles incident to points of the sector. By Theorem \ref{boundedsidescircle} these circles must be incident to a pair of points in a cell.

\qed

\begin{rmk}\label{rmk:toycase} Figure \ref{fig:twobushcircles} contains claims about the points and the centers of the unit circles being contained in rectangles of width $\delta$ and length $\sqrt{\delta}$ which are distance 1 away and parallel along their long axis. These claims are true. 
\end{rmk}

\noindent \textit{Proof.} By Lemma \ref{lem:circlesNminusball} we can restrict the point set and center point of the unit circles to each lie in some ball of radius $\sqrt{\delta} \lesssim N^-$ at distance 1 from each other. Let's call these balls $B_1$ and $B_2$.
    
By Lemma \ref{lem:circle_arc_tube_sqrt} we know that the circle arcs from the two organizing circles which intersect $B_2$ are each contained in some tube of width $\delta$. The corresponding bush circles $\{c_{1,i}\}$ have centers $\lesssim \sqrt{\delta}$ away from each other so the angles of the tubes which contain each of the bush circle arcs in $B_1$ differ by at most $\sqrt{\delta}$. Also these tubes are all concurrent at $p_1$. Same for the bush circles $\{c_{2,i}\}$ which are all concurrent at $p_2$. The tubes have width $\delta$ and length $\sqrt{\delta}$ so a segment of length $\sqrt{\delta}$ can have angle as much as $\delta/\sqrt{\delta} = \sqrt{\delta}$ away from the axis of the tube and still be entirely contained in the tube. Thus the bush circle arcs $\{c_{1,i}\}$, respectively $\{c_{1,i}\}$ in $B_1$ are contained in a single tube $T_1$, respectively $T_2$ of width $\delta$ and length $\sqrt{\delta}$. The cells of the two bush cell decomposition from Theorem \ref{doublebushmixingcircle} are contained in $T_1 \cap T_2$ so they are contained in some $\delta \times \sqrt{\delta}$ tube.

\qed

\subsection{A proto inverse theorem for unit distances} \label{section:protoinverseSTcircle}

A proto inverse Szemer\'edi-Trotter theorem for unit circles will be a recipe for constructing a configuration
of points and unit circles which may not terminate or may not yield an extremal configuration but
so that a large portion of every extremal example can be obtained using this recipe.

When considering such a recipe, an important piece of information is how many parameters
one needs to specify to obtain an instantiation of the recipe.

There is a trivial recipe taking $O(N)$ parameters. Namely use $4N$ parameters to 
completely specify a set of $N$ unit circles, $\lines$ and a set of $N$ points $\pts$. Examine
the set of incidences between these unit circles and points $I(\lines,\pts)$. If it happens
to be that $|I(\lines,\pts)| \gtrsim N^{{4 \over 3}-}$, then we have constructed
an extremal configuration and in fact, every extremal configuration can be constructed in 
this way. This recipe and its proto inverse Szemer\'edi-Trotter theorem amount to really
just the definition of extremal configuration, and nothing has been gained.

Now, however, we describe a recipe using just $O(N^{{1 \over 3}})$ parameters.
Our recipe will be based on the two bush cell decomposition obtained from Theorem \ref{doublebushmixingcircle}.  $a_1 < a_2  < \dots < a_{N^{{1 \over 3}}}$ will be real numbers
representing the slopes of the lines tangent to the circles $\{c_{1,j}\}$ from the first bush at the first organizing point $p_1$. $b_1 < b_2  < \dots < b_{N^{{1 \over 3}}}$
will be real numbers
representing the slopes of the lines tangent to the circles $\{c_{2,j}\}$ from the second bush at the second organizing point $p_2$. The final ingredients will be a set of unit circles
$c_{s,1} \dots  c_{s,N^{{1 \over 3}}}$ which will serve as the structuring circles for the
cells in the sector bounded by $c_{1,1}$ and $c_{1,2}$. We denote by $\cell_{i,j}$ the cell that is the intersection of the sector bounded by $c_{1,i}$ and $c_{1,i+1}$ and of the sector bounded by $c_{2,j}$ and $c_{2,j+1}$.

We declare the recipe to have failed unless there are at least $N^{{1 \over 3}-}$ values
of $j$ so that at least $N^{{ 1 \over 3}-}$ and no more than $N^{{1 \over 3}+}$ of the
crossings between the structuring circles $c_s$ are in the cell $\cell_{1,j}$. If the recipe did not fail yet, we define these crossings to be the points of the
cell $\cell_{1,j}$. For each cell, we find all unit circles which are incident to
two points. We declare this set of circles to be $\lines$. We say that the recipe has failed
unless there are at least $N^{{2 \over 3}-}$ choices of $(i,j)$ so that $N^{{2 \over 3} \pm}$ circles of $\lines$ cross the cell $\cell_{i,j}$. Otherwise, we
say that the recipe has failed unless for at least $N^{{2 \over 3}-}$ of these choices
$(i,j)$ the circles going through the $(i,j)$th cell are structured. If they are structured,
we refer to the structuring points as the points of the $(i,j)$th cell. And combining all
of these structuring points, we get the set $\pts$ and we declare that the construction
has succeeded. In this case, $(\lines,\pts)$ is an extremal configuration. We denote the
output of this recipe as a function of its inputs as $(\lines(a,b,l_s), \pts(a,b,l_s))$.

\begin{thm}  \label{protoinverseSTcircles}
Let $(\lines,\pts)$ be an extremal unit circle configuration. Then there is a choice of the 
$O(N^{{1 \over 3}})$ parameters $(a,b,l_s)$ so that
$(P(\lines) \cap \lines(a,b,l_s),P(\pts) \cap \pts(a,b,l_s))$ is an extremal configuration
(of $N^{1-}$ unit circles $N^{1-}$ points and $N^{{4 \over 3}-}$ incidences.)
\end{thm}

\begin{proof}  This is essentially a consequence of Theorem \ref{doublebushmixingcircle} and
its proof. We find
a sector through $p_1$ through which $N^{1-}$ circles make $N^{1-}$ incidences with
the structured points total and each circle making at least two incidences per cell. (Possibly by removing all
but $N^{-}$ density of elements of each bush.) Then we choose $a,b$ so that they give
all the cells of the double bush construction with $a_1,a_2$ corresponding to the sector chosen above and we let $l_s$ be the structuring circles for the sector. \end{proof}

\bigskip

Recall that an organizing circle is partitioned into $N^{{1 \over 3}\pm}$ intervals by points it is incident to and a $N^+$ density of the intervals intersect a structured set of $N^{{2 \over 3}\pm}$ circles. These sets of structured circles are dual to the sets of points contained in each sector defined by an organizing point.

\begin{defn} Let $(\lines,\pts)$ be an extremal configuration. Given an organizing circle $c$ and its bush of points from $\pts$: $\{p_{1}, \dots, p_{M_1}\}$ on $c$ with $M\gtrsim N^{{1 \over 3}-}$, for each $j \in [1, M]$ we call the subset of $\lines$ which intersects $l_1$ between $p_{j,1}$ and $p_{j+1,1}$ a sector of circles. Note each organizing point defines $\gtrsim N^{{1 \over 3}-}$ many sectors of circles.
\end{defn}

In Theorem \ref{doublebushmixingcircle} we saw that two organizing points and their bush of circles yield a point-weighted cell decomposition where the cells are the intersection of a sector from the first organizing point with a sector from the second organizing point. We similarly define the dual notion of cells of unit circles:

\begin{defn}Let $(\lines,\pts)$ be an extremal configuration. Given a pair of organizing circles $c_1$ and $c_2$ in $\lines$ and their bush of points from $\pts$: $\{p_{1,1}, \dots, p_{M_1,1}\}$ on $c_1$ and $\{p_{1,2}, \dots, p_{M_2,2}\}$ on $c_2$ with $M_1,M_2 \gtrsim N^{{1 \over 3}-}$ we define $\gtrsim N^{{2 \over 3}-}$ many cells of circles as follows. For each $(j,k) \in [1, M_1] \times [1, M_2]$ the cell of circles $(j,k)$ is the set of circles in $\lines$ that intersects $l_1$ in the interval between $p_{j,1}$ and $p_{j+1,1}$ and intersects $l_2$ in the interval between $p_{k,2}$ and $p_{k+1,2}$.
\end{defn}

\begin{thm} \label{dualprotoinverseSTcircle}
Let $(\lines,\pts)$ be an extremal unit circle configuration. Then there are $\gtrsim N^{2-}$ pairs of organizing circles $c_1$, $c_2$ in $\lines$ and a refinement $\lines' \subset \lines$ so that $(c_1, c_2)$ partition $\lines'$ into $\gtrsim N^{{2 \over 3}-}$ many cells of unit circles each of cardinality $\sim N^{{1 \over 3}\pm}$. Furthermore each of the sectors of circles are structured by $N^{{1 \over 3} \pm}$ structuring points.
\end{thm}

\begin{proof}  Apply point-unit circle duality to Theorem \ref{doublebushmixingcircle} (a unit circle is mapped to its center point and a point to the unit circle centered at that point) with Lemma \ref{circleduality}.
\end{proof}

\begin{thm}[Mixing]
\label{mixingcircle}
Let $(\lines,\pts)$ be an extremal unit circle configuration with the cell decomposition given by Theorem \ref{protoinverseSTcircles}. Then $\gtrsim N^{{4 \over 3}-}$ pairs of cells share $\gtrsim N^{{1 \over 3}-}$ circles which take at least two incidences of $J'$ in each of the two cells.
\end{thm}

\begin{proof}
    Assume we have the cell decomposition from Theorem \ref{protoinverseSTcircles}. An application of Theorem \ref{dualprotoinverseSTcircle} gives us that $\gtrsim N^{2-}$ pairs of circles are organizing. 

    Consider a pair of organizing circles $c_1, c_2$. We take a refinement of the circle set such that each circle has at least two incidences in a cell that $c_1$ has incidences in and at least two incidences in a cell that $c_2$ has incidences in. Note this must yield an extremal refinement for $\gtrsim N^{2-}$ pairs of organizing circles $(c_1,c_2)$ because each of the $\gtrsim N^{1-}$ organizing circles in $\lines'$ has at least 2 incidences in $\gtrsim N^{{1 \over 3}-}$ cells which each have $\gtrsim N^{{2 \over 3}-}$ circles going through it which each have at least two incidences in that cell. So each organizing circle contributes $\gtrsim N^{1-}$ circles (out of a total of $< N$ circles) so $\gtrsim N^{1-}$ other organizing circles must share $\gtrsim N^{1-}$ regular circles that take at least two points in one of each of their cells.

    From our initial application of Theorem \ref{dualprotoinverseSTcircle} we know that $\gtrsim N^{{2 \over 3}-}$ pairs of circle arcs between adjacent points on $c_1$ and adjacent points on $c_2$ share $\gtrsim N^{{1 \over 3}-}$ circles. (Note this is still true after our refinement because we kept $\gtrsim N^{1-}$ circles). Furthermore circles take $\lesssim N^{+}$ incidences in a single cell so $\gtrsim N^{{2 \over 3}-}$ pairs of cells that $c_1$ and $c_2$ take incidences in share $\gtrsim N^{{1 \over 3}-}$ circles out of $\lesssim N^{{2 \over 3}}$ total pairs of cells. 

    Since this result holds for any generic pair of $\gtrsim N^{2-}$ pairs of organizing circles, we conclude that any generic pair of cells must share $\gtrsim N^{{1 \over 3}-}$ circles.
\end{proof}

\end{document}